\documentclass[leqno,12pt]{article}

\usepackage{amssymb}
\usepackage{euscript}
\usepackage{natbib}

\evensidemargin\oddsidemargin

\begin{document}

\baselineskip=18pt \setcounter{page}{1}

\renewcommand{\theequation}{\thesection.\arabic{equation}}
\newtheorem{theorem}{Theorem}[section]
\newtheorem{lemma}[theorem]{Lemma}
\newtheorem{proposition}[theorem]{Proposition}
\newtheorem{corollary}[theorem]{Corollary}
\newtheorem{remark}[theorem]{Remark}
\newtheorem{fact}[theorem]{Fact}
\newtheorem{problem}[theorem]{Problem}
\newtheorem{corollaire}[theorem]{Corollaire}
\newcommand{\eqnsection}{
\renewcommand{\theequation}{\thesection.\arabic{equation}}
    \makeatletter
    \csname  @addtoreset\endcsname{equation}{section}
    \makeatother}
\eqnsection

\def\r{{\mathbb R}}
\def\s{{\mathcal S}}
\def\e{{\mathbb E}}
\def\p{{\mathbb P}}
\def\bZ{{\mathbb Z}}
\def\S{{\mathbb S}}
\def\t{{\mathbb T}}
\def\bN{{\mathbb N}}
\def\deg{{\rm b}}
\def\P{{\bf P}}
\def\ind {\mbox{\rm 1\hspace {-0.04 in}I}}
\def\E{{\bf E}}
\def\ee{\mathrm{e}}
\def\d{\, \mathrm{d}}
\newcommand{\cF}{\mathcal{F}}
\newcommand{\cG}{\mathcal{G}}
\newcommand{\cU}{\mathcal{U}}
\def\Q{{\bf Q}}
\def\q{{\mathbb Q}}
\def\bp{{\bf p}}
\def\tS{{\widetilde{S}}}
\def\tQ{{\widetilde{Q}}}
\def\tX{{\widetilde{X}}}
\newcommand{\cT}{\mathcal{T}}

\def\Z{\mathbb{Z}}
\def\p{\mathbb{P}}
\def\V{\mathbb{V}}
\def\t{\mathbb{T}}
\def\ind {\mbox{\rm 1\hspace {-0.04 in}I}}


\vglue15pt

\centerline{\Large\bf Large deviations for transient random walks}

\bigskip

\centerline{\Large\bf in random environment on a Galton--Watson
tree}
\bigskip
\bigskip
\bigskip
\centerline{by}

\bigskip

\centerline{Elie Aid\'ekon}
\medskip

\centerline{\it Universit\'e Paris VI}

\bigskip
\bigskip
\bigskip
\bigskip

{\leftskip=2truecm \rightskip=2truecm \baselineskip=15pt \small

\noindent{\slshape\bfseries Summary.} Consider a random walk in
random environment on a supercritical Galton--Watson tree, and let
$\tau_n$ be the hitting time of generation $n$. The paper presents a
large deviation principle for $\tau_n/n$, both in quenched and
annealed cases. Then we investigate the subexponential situation,
revealing a polynomial regime similar to the one encountered in one
dimension. The paper heavily relies on estimates on the tail
distribution of the first regeneration time.

\bigskip


\bigskip

} 

\bigskip
\bigskip


\section{Introduction}

We consider a super-critical Galton--Watson tree $\t$  of root $e$
and offspring distribution $(q_k,\,k\ge 0)$ with finite mean  $
m:=\sum_{k\ge 0}kq_k>1$. For any vertex $x$ of $\t$, we call $|x|$
the generation of $x$, $(|e|=0)$ and $\nu(x)$ the number of children
of $x$; we denote these children by $x_i,\,1\le i\le \nu(x)$. We let
$\nu_{min}$ be the minimal integer such that $q_{\nu_{min}}>0$ and
we suppose that $\nu_{min}\ge 1$ (thus $q_0=0)$. In particular, the
tree survives almost surely. Following Pemantle and Peres
\cite{PemPer}, on each vertex $x$, we pick independently and with
the same distribution a random variable $A(x)$, and we define
\begin{itemize}
\item $\omega(x,x_i):={A(x_i)\over 1+\sum_{i=1}^{\nu(x)}A(x_i)}$\,, $\forall\, 1\le i\le \nu(x)$,
\item $\omega(x,\buildrel \leftarrow  \over x):={1\over 1+\sum_{i=1}^b A(x_i)}$\,.
\end{itemize}

\noindent To deal with the case $x=e$, we add a parent ${\buildrel
\leftarrow \over e}$ to the root and we set $\omega({\buildrel
\leftarrow \over e},e)=1$. Once the environment built, we define the
random walk $(X_n,\,n\ge 0)$ starting from $y\in \t$ by
\begin{eqnarray*}
P_{\omega}^y(X_0=y) &=&1,\\
P_{\omega}^y(X_{n+1}=z \, | \, X_n=x)&=&\omega(x,z)\,.
\end{eqnarray*}

\noindent The walk $(X_n,\,n\ge 0)$ is a $\t$-valued Random Walk in
Random Environment (RWRE). To determine the transience or recurrence
of
the random walk,  Lyons and Pemantle \cite{lp92} provides us with the following criterion. Let $A$ be a generic random variable having the distribution of $A(e)$.\\

\noindent {\bf Theorem A (Lyons and Pemantle~\cite{lp92})} {\it The walk $(X_n)$ is transient if $\,\inf_{[0,1]} \E[A^t]> \frac{1}{m}$, and is recurrent otherwise.}\\

In the transient case, let $v$ denote the speed of the walk, which
is the deterministic real $v\ge 0$ such that
$$
\lim_{n\rightarrow \infty} \frac{|X_n|}{n}=v, \qquad a.s.
$$

\noindent Define
\begin{eqnarray*}
i  &:=&   \mbox{ess inf}\, A \,,\\
s  &:=&   \mbox{ess sup} \,A\,.
\end{eqnarray*}

\noindent We make the hypothesis that $0<i \le s<\infty$. Under this
assumption, we gave a criterion in \cite{aid07} for the positivity
of the speed $v$. Let

\begin{equation}
    \Lambda:=Leb \left\{t\in \r: \E[A^t]\le
    {1\over q_1}\right\} \qquad (\Lambda=\infty \; \mbox{if}\; q_1=0).
    \label{LD:Lambda}
\end{equation}

\noindent {\bf Theorem B~ (\cite{aid07}) }{\it
 Assume $\inf_{[0,1]} \E[A^t]> \frac{1}{m}$, and let
 $\Lambda$ be as in $(\ref{LD:Lambda})$.

 {\rm (a)} If $\Lambda< 1$, the walk has zero speed.

 {\rm (b)} If $\Lambda>1$, the walk has positive
 speed. }\\

When the speed is positive, we would like to have information on
how hard it is for the walk to have atypical behaviours, which means
to go a little faster or slower than its natural pace. Such
questions have been discussed in the setting of biased random walks
on Galton--Watson trees, by Dembo et al. in \cite{dgpz02}. The
authors exhibit a large deviation principle both in quenched and
annealed cases. Besides, an uncertainty principle allows them to
obtain the equality of the two rate functions. For the RWRE in
dimensions one or more, we refer to Zeitouni \cite{zeitouni04} for a
review of the subject. In our case, we consider a random walk which
always avoids the parent ${\buildrel \leftarrow \over e}$ of the
root,
 and we obtain a large deviation principle, which, following
\cite{dgpz02}, has been divided into two parts.

We suppose in the rest of the paper that
\begin{eqnarray}
\inf_{[0,1]} \E[A^t] &>& \frac{1}{m}\,, \label{LD:trans}\\
\Lambda &>& 1 \label{LD:posispeed}\,,
\end{eqnarray}

\noindent which ensures that the walk is transient with positive
speed. Before the statement of the results, let us introduce some
notation. Define for any $n\ge 0$ and $x\in \t$,
\begin{eqnarray*}
   \tau_n      &:=&       \inf\left\{k\ge 0\, : |X_k|=n\right\}\,, \\
    D(x)
 &:=&
            \inf\left\{k\ge 1\,: X_{k-1}=x,\, X_k={\buildrel \leftarrow
                     \over x} \right\},~~(\inf\emptyset:=\infty)\,.
\end{eqnarray*}

\noindent Let $\P$ denote the distribution of the environment
$\omega$ conditionally on $\t$, and $\Q:=\int \P(\cdot) GW(\! \d \t)$. Similarly, we denote by $\p^x$ the distribution defined by $\p^x
(\cdot) := \int P_\omega^x (\cdot) \P(\! \d \omega)$ and
by $\q^x$ the distribution
\begin{eqnarray*}
  \q^x(\cdot) := \int \p^x(\cdot) GW(\! \d \t)\,.
\end{eqnarray*}

\begin{theorem}
\label{LD:speedup} {\bf (Speed-up case)} There exist two continuous,
convex and strictly decreasing functions $I_a \le I_q$ from
$[1,1/v]$ to $\mathbb{R}_+$, such that $I_a(1/v)=I_q(1/v)=0$ and for
$ a<b$, $b\in [1, 1/v]$,
\begin{eqnarray}
     \lim_{n \rightarrow \infty} { 1 \over n }\ln \q^e\left({\tau_n\over n}\in ]a, b] \right)
&=&
     - I_a(b)\,, \label{LD:devupaup}\\
     \lim_{n\rightarrow\infty} {1\over n}\ln P_{\omega}^e\left( {\tau_n \over n} \in ]a, b] \right)
&=&
     - I_q(b)\,. \label{LD:devupqup}
\end{eqnarray}
\end{theorem}

\begin{theorem}
\label{LD:slowdown} {\bf (Slowdown case)} There exist two
continuous, convex functions $I_a \le I_q$ from $[1/v,+\infty[$ to
$\mathbb{R}_+$, such that $I_a(1/v)=I_q(1/v)=0$ and for any $1/v \le
a < b$,
\begin{eqnarray}
     \lim_{n\rightarrow\infty} {1\over n}\ln \q^e\left({\tau_n \over n} \in [a,b[\right)
&=&
     -I_a(a)\,, \label{LD:devupaslow}\\
     \lim_{n\rightarrow\infty} {1\over n}\ln P_{\omega}^e\left({\tau_n \over n}\in [a,b[\right)
&=&
     -I_q(a)\,. \label{LD:devupqslow}
\end{eqnarray}
Besides, if $i> \nu_{min}^{-1}$, then $I_a$ and $I_q$ are strictly
increasing on $[1/v,+\infty[$. When $i\le \nu_{min}^{-1}$, we have
$I_a=I_q=0$ on the interval.
\end{theorem}

\noindent As pointed by an anonymous referee, it would be
interesting to know when $I_a$ and $I_q$ coincide. We do not know
the answer in general. However, the computation of the value of the
rate functions at $b=1$ reveals situations where the rate functions
differ. Let
\begin{eqnarray*}
    \psi(\theta)
:=
    \ln\left( E_{\Q}\left[\sum_{i=1}^{\nu(e)} \omega(e,e_i)^{\theta} \right]\right). \label{LD:psi}
\end{eqnarray*}

\noindent Then $\psi(0)= \ln(m) $ and $\psi(1)=\ln\left(E_{\Q}\left[
\sum_{i=1}^{\nu(e)} \omega(e,e_i) \right] \right)$.

\begin{proposition}
\label{LD:i1} We have
\begin{eqnarray}
    I_a(1)
&=&
    -\psi\left(1\right)\,, \label{LD:ia1}\\
    I_q(1)
&=&
    -\inf_{]0,1]}{1\over \theta} \psi(\theta)\,. \label{LD:iq1}
\end{eqnarray}
In particular, $I_a(1)=I_q(1)$ if and only if $\psi'(1)\le \psi(1)$.
Otherwise $I_a(1)<I_q(1)$.
\end{proposition}

Quite surprisingly, we can exhibit elliptic environments on a
regular tree for which the rate functions differ. This could hint that
the  uncertainty of the location of the first passage in \cite{dgpz02} does
not hold anymore for a random environment. Here is an explicit example. Consider a binary tree ($q_2=1$). Let $A$ equal $0.01$ with probability $0.8$ and $500$ with probability $0.2$. Then we check that the walk is transient, but $\psi'(1)>\psi(1)$ so that $I_a(1)\neq I_q(1)$ on such an environment.

\bigskip

Theorem \ref{LD:slowdown} exhibits a subexponential regime in the
slowdown case when $i\le \nu_{min}^{-1}$. The following theorem
details this regime. Let
\begin{eqnarray*}
            \S^e(\cdot)
&:=&
            \q^e(.\,|\, D(e)=\infty)\,.
\end{eqnarray*}

\begin{theorem}
\label{LD:subdev}
We place ourself in the case $i< \nu_{min}^{- 1} $. \\
(i) Suppose that either ``$i< \nu_{min}^{- 1} $ and  $q_1=0$'' or
``$i< \nu_{min}^{- 1} $ and  $s < 1$" . There exist constants
$d_1,d_2 \in (0,1)$ such that for any $a>1/v$ and $n$ large enough,
\begin{eqnarray}
\label{LD:subeq} e^{-n^{d_1}} < \S^e( \tau_n >a n)<e^{-n^{d_2}}\,.
\end{eqnarray}
\\
(ii) If $q_1>0$ and $s>1$ (id est when $\Lambda<\infty$), the regime
is polynomial and we have for any $a>1/v$,
\begin{eqnarray}
\label{LD:poleq} \lim_{n\rightarrow\infty } {1\over
\ln(n)}\ln\left(\S^e( \tau_n
>a n)\right)= 1-\Lambda\,.
\end{eqnarray}
\end{theorem}

We mention that in one dimension, which can be seen as a critical
state of our model where $q_1=1$, such a polynomial regime is proved
by Dembo et al. \cite{dpz96}, our parameter $\Lambda$ taking the
place of the
well-known $\kappa$ of Kesten, Kozlov, Spitzer \cite{Kes75}. We did not deal with the critical case $i=\nu_{min}^{-1}$. Furthermore, we do not have any conjecture on the optimal values of $d_1$ and $d_2$ and do not know if the two values are equal.\\

The rest of the paper is organized as follows. Section 2 describes
the tail distribution of the first regeneration time, which is a
preparatory step for the proof of the different theorems. Then we
prove Theorems \ref{LD:speedup} and \ref{LD:slowdown} in Section 3,
which includes also the computation of the rate functions at speed
$1$ presented in Proposition \ref{LD:i1}. Section 4 is devoted to
the subexponential regime with the proof of Theorem \ref{LD:subdev}.


\section{Moments of the first regeneration time}

We define the first regeneration time
$$
\Gamma_1:=\inf\left\{k>0 \, : \nu(X_k)\ge 2,\,
D(X_k)=\infty,\,k=\tau_{|X_k|} \right\}
$$

\noindent as the first time when the walk reaches a generation by a
vertex having more than two children and never returns to its
parent. We propose in this section to give information on the tail
distribution  of $\Gamma_1$ under $\S^e$. We first introduce some
notation used throughout the paper. For any $x\in\t$, let
\begin{eqnarray}
N(x)       &:=&       \sum_{k\ge 0}\ind_{\{X_k=x\}}\,, \label{LD:green}\\
T_x         &:=&       \inf\left\{k\ge 0 :
                                           \; X_k=x \right\}\,, \nonumber \\
T_x^*     &:=&  \inf\{k\ge 1 \,:\, X_k=x\} \nonumber \,.
\end{eqnarray}

\noindent This permits to define
\begin{eqnarray}
\beta(x)  &:=&  P_{\omega}^x(T_{{\buildrel \leftarrow \over
x}}=\infty)\,,\nonumber \\
\gamma(x) &:=&  P_{\omega}^x(T_{{\buildrel \leftarrow \over
x}}=T_x^*=\infty)\,. \label{LD:gamma}
\end{eqnarray}

\noindent The following fact can be found in \cite{dgpz02} (Lemma 4.2) in the case of biased random walks, and is directly adaptable in our setting.\\ \\
{\bf Fact A} The first regeneration height $|X_{\Gamma_1}|$ admits
exponential
moments under the measure $\S^e(\cdot)$.\\

\subsection{The case $i > \nu_{min}^{-1}$}

This section is devoted to the case $i > \nu_{min}^{-1}$, where
$\Gamma_1$ is proved to have exponential moments.

\begin{proposition}
\label{LD:timexp} Suppose that $i>\nu_{min}^{-1}$. There exists
$\theta>0$ such that $E_{\S^e}\left[e^{\theta
{\Gamma_1}}\right]<\infty$.
\end{proposition}
\noindent {\it Proof}. The proof follows the strategy of Proposition
1 of Piau \cite{Piau98}. We couple the distance of our RWRE to the
root $(|X_n|)_{n\ge 0}$ with a biased random walk $(Y_n)_{n\ge 0}$
on $\bZ$ as follows. Let $ p : = {i \nu_{min} \over 1 + i \nu_{min}
} $, and let $u_n,\,n\ge 0$, be a family of i.i.d.$\!$ uniformly
distributed [0,1] random variables. We set $X_0=e$ and $Y_0=0$. If
$X_k$ and $Y_k$ are known, we construct
\begin{eqnarray*}
X_{k+1}     &=&   x_i    \qquad \qquad \qquad \qquad   \mbox{if}\;\; \sum_{\ell=1}^{i-1}\omega(x,x_{\ell})\le u_k<\sum_{\ell=1}^{i}\omega(x,x_{\ell})\,,\\
X_{k+1}     &=&     {\buildrel \leftarrow \over x}     \qquad \qquad \qquad \qquad             \mbox{otherwise}\,,\\
Y_{k+1}     &=&        y+2\ind_{\{u_k\le p\}}-1\,,
\end{eqnarray*}

\noindent where $ x : = X_k \in \t $ and $ y : = Y_k \in \bZ $. Then
$(X_n)_{n\ge 0}$ has the distribution of our $\t$-RWRE indeed, and
$(Y_n)_{n\ge 0}$ is a random walk on $\bZ$ which increases of one
unit with probability $p>1/2$ and decreases of the same value with
probability $1-p$. Notice also that on the event $\{D(e)=\infty\}$,
we have
$$
|X_{k+1}|-|X_k|  \ge   Y_{k+1}-Y_k\,.
$$

\noindent It implies that the first regeneration time
$\mathcal{R}_1$ of $(Y_n)_{n\ge 0}$ defined by
$$
\mathcal{R}_1:=\inf\left\{k>0 \, : Y_{\ell}< Y_k \;\forall \ell < k\,, Y_{m}\ge Y_k\; \forall m>k
 \right\}
$$
is necessarily a regeneration
time for $(X_n,\,n\ge 0)$, which proves in turn that
\begin{eqnarray*}
\S^e(\Gamma_1>n) \le \q^e(\mathcal{R}_1>n)\,.
\end{eqnarray*}

\noindent To complete the proof, we must ensure that
$\q^e(\mathcal{R}_1>n)$ is exponentially small, which is done in
\cite{dpz96} Lemma 5.1. $\Box$

\subsection{The cases ``$i<\nu_{min}^{-1}$, $q_1=0$" and `` $i<\nu_{min}^{-1}$, $s<1$"}

When $i<\nu_{min}^{-1}$, if we assume also that $q_1=0$ or $s<1$, we
prove that $\Gamma_1$ has a subexponential tail. This situation
covers, in particular, the case of RWRE on a regular tree.

\begin{proposition}
\label{LD:timereg} Suppose that $i<\nu_{min}^{-1}$ and $ q_1=0 $,
then there exist $1 > \alpha_1 > \alpha_2 >0$ such that for $n$
large enough,
\begin{equation}
e^{-n^{\alpha_1}}<\S^e(\Gamma_1>n)<e^{-n^{\alpha_2}}\,.
\label{LD:sub}
\end{equation}
The same relation holds with some $1>\alpha_3>\alpha_4>0$ in the
case ``$i<\nu_{min}^{-1}$ and $s<1$".
\end{proposition}
\noindent {\it Proof of Proposition \ref{LD:timereg}: lower bound}.
We only suppose that $i<\nu_{min}^{-1}$, which allows us to deal
with both cases of the lemma. Define for some $p' \in (0,1/2)$ and
$b\in \bN$,
\begin{eqnarray*}
w_+  &:=&    \Q\left(\sum_{i=1}^{\nu}A(e_i)\ge {1-p'\over p'},\,\nu(e) \le b\right)\,, \\
w_-  &:=&    \Q\left(\sum_{i=1}^{\nu}A(e_i)\le {p'\over 1-p'}, \, \nu(e) \le
b\right)\,.
\end{eqnarray*}

\noindent By (\ref{LD:trans}),
$E_{\Q}\left[\sum_{i=1}^{\nu(e)}A(e_i) \right]>1$ and therefore
$\Q\left(\sum_{i=1}^{\nu(e)}A(e_i)>1\right)>0$. Since $\mbox{ess inf
}A<\nu_{min}^{-1}$, it guarantees that
$\Q\left(\sum_{i=1}^{\nu(e)}A(e_i) <1\right)>0$. Consequently, by
choosing $p'$ close enough of $1/2$ and $b$ large, we can take $w_+$
and $w_-$ positive.
 Let $c:={1\over 6 \ln(b)}$, and define $h_n:=\lfloor c\ln(n)\rfloor$. A tree $\t$ is said to be $n$-good if
\begin{itemize}
\item any vertex $x$ of the $h_n$ first generations verifies
 $\nu(x)\le b$ and $\sum_{i=1}^{\nu(x)}A(x_i)\ge {1-p'\over p'}$\,,
\item any vertex $x$ of the $h_n$ following generations verifies $\nu(x)\le
b$ and $\sum_{i=1}^{\nu(x)}A(x_i)\le {p'\over 1-p'}$\,.
\end{itemize}
\bigskip
\noindent We observe that $\Q(\t \mbox{ is}\; n\mbox{-good}) \ge
w_+^{h_n b^{h_n}} w_-^{h_n b^{2h_n}}\ge e^{-n^{1/3+o(1)}}$ which is
stretched exponential, i.e. behaving like $e^{-n^{r+o(1)}}$ for some
$r\in (0,1)$. Define the events
\begin{eqnarray*}
E_1     &:=&            \{ \mbox{at time }  \;  \tau_{h_n}
\;\mbox{we can't find an edge of
                          level smaller than}
                                           \;   h_n   \;  \mbox{crossed only once}\}\\
         &&\; \cap \;\{
                                                                D(e)>\tau_{h_n}
                            \}\,,\\
E_2     &:=&            \{\mbox{the walk visits the level} \; h_n
                           \; n \;\mbox{times before reaching the root or the level}\; 2h_n
                           \}\,,\\
E_3     &:=&            \{\mbox{after  the}\;n\mbox{-th visit of
                           level}\;h_n, \;\mbox{the walk reaches level} \;2h_n \; \mbox{before level} \; h_n
                           \}\,, \\
E_4     &:=&            \{ \mbox{after time}\; \tau_{2h_n}\;
\mbox{the
                            walk never
                           comes back to level}\; 2h_n-1 \}\,.
\end{eqnarray*}

\noindent Suppose that the tree is $n$-good. Since $A$ is supposed
bounded, there exists a constant $c_1>0$ such that for any $x$
neighbour of $y$, we have
\begin{equation}
\label{LD:omega} \omega(x,y)\ge {c_1\over \nu(x)}\,.
\end{equation}

\noindent It yields that $P_{\omega}^e(E_1)^{-1}=O(n^K)$ for some
$K>0$ (where $O(n^K)$ means that the function is bounded by a factor
of $n\rightarrow n^K$). Combine (\ref{LD:omega}) with the strong
Markov property at time $\tau_{h_n}$ to see that
\begin{eqnarray*}
P_{\omega}^e(E_3\,|\, E_1 \cap E_2)^{-1}=O(n^K)\,,
\end{eqnarray*}

\noindent where $K$ is taken large enough. We emphasize that the
functions $O(n^K)$ are deterministic. Still by Markov property,
\begin{eqnarray}
P_{\omega}^e(E_1\cap E_2 \cap E_3 \cap E_4) =
E_{\omega}^e[\ind_{E_1\cap E_2 \cap E_3}\beta(X_{\tau_{2h_n}})]\,.
\label{LD:e1234}
\end{eqnarray}

\noindent Let $(Y_n')_{n\ge 0}$ be the random walk on $\bZ$ starting
from zero with
$$P_{\omega}(Y_{n+1}'=k+1\,|\,Y_{n}'=k)=1-P_{\omega}(Y_{n+1}'=k-1\,|\,Y_{n}'=k)=p'\,.$$

\noindent We introduce $T_i':=\inf\{k\ge 0 \,: \, Y_k=i\}$, and
$p_n'$ the probability that $(Y_n')_{n\ge 0}$ visits $h_n$ before
$-1$:
$$
p_n':=P_{\omega}(T_{-1}'< T_{h_n}')\,.
$$

\noindent By a coupling argument similar to that encountered in the
proof of Proposition \ref{LD:timexp}, we show that in an $n$-good
tree,
\begin{eqnarray}
\label{LD:nK} P_{\omega}^e(E_1 \cap E_2) \;\ge \;
P_{\omega}^e(E_1)(p_n')^n \; =  \;  O(n^K)^{-1}(p_n')^n \,,
\end{eqnarray}

\noindent which gives
\begin{eqnarray}
P_{\omega}^e(E_1 \cap E_2 \cap E_3) \ge O(n^K)^{-1}(p_n')^n\,.
\label{LD:e123'}
\end{eqnarray}

\noindent Observing that $\q^e(\Gamma_1>n,\, D(e)=\infty) \ge
E_{\Q}\left[\ind_{\{\t \;\mbox{is} \; n\mbox{-good}\}}\ind_{E_1\cap
E_2\cap E_3 \cap E_4}\right]$, we obtain by (\ref{LD:e1234})
\begin{eqnarray*}
       \q^e(\Gamma_1>n,\,D(e)=\infty)
&\ge&
       E_{\q^e}\left[\ind_{\{\t \;\mbox{is n-good}\}}\ind_{E_1 \cap E_2 \cap E_3}\beta(X_{\tau_{2h_n}})\right]\\
&=&
       E_{\q^e}\left[\ind_{\{\t \;\mbox{is n-good}\}}P_{\omega}^e(E_1 \cap E_2 \cap
       E_3)\right]E_{\Q}[\beta]\,,
\end{eqnarray*}

\noindent by independence. By (\ref{LD:e123'}),
\begin{eqnarray*}
      \q^e(\Gamma_1>n,\,D(e)=\infty)
\ge
      O(n^K)^{-1}\Q\left(\t \;\mbox{is} \;
       n\mbox{-good}\right)(p'_n)^n\,.
\end{eqnarray*}

\noindent We already know that $\Q\left(\t \;\mbox{is} \;
       n\mbox{-good}\right)$ has a stretched exponential lower bound, and it remains to observe that the same holds for $(p'_n)^n$. But
the method of gambler's ruin shows that $p_n'\ge 1-\left({p'\over
1-p'}\right)^{h_n}$, which gives the required lower
bound by our choice of $h_n$. $\Box$\\

Let us turn to the upper bound. We divide the proof in two,
depending on which case we deal with.\\

\noindent {\it Proof of Proposition \ref{LD:timereg}: upper bound in
the case $q_1=0$}. Assume that $q_1=0$ (the condition
$i<\nu_{min}^{-1}$ is not required in the proof). The proof of the
following lemma is deferred. Recall the notation introduced in
(\ref{LD:gamma}), $\gamma(e):=P_{\omega}^e(T_{{\buildrel \leftarrow
\over e}}=T_e^*=\infty)\le \beta(e)$.
\begin{lemma}
\label{LD:beta} When $q_1=0$, there exists a constant $c_2 \in(0,1)$
such that for large $n$, $$
E_{\Q}\left[\left(1-\gamma(e)\right)^n\right]\le e^{-n^{c_2}} \,. $$
\end{lemma}

\noindent Denote by $\pi_k$ the $k$-th distinct
site visited by the walk $(X_n,\, n\ge 0)$. We observe that
\begin{eqnarray}
       \q^e(\Gamma_1>n^3)
&\le&
       \q^e(\Gamma_1>\tau_{n}) + \q^e(\, \mbox{more than } \;
       n^2 \;   \mbox{distinct sites are visited before}\; \tau_n) \nonumber \\
&&
    + \,\, \q^e(\exists \, k\le n^2\,:\,N(\pi_k)> n). \label{LD:regen1}
\end{eqnarray}

\noindent Since $\q^e(\Gamma_1 > \tau_{n}) = \q^e(|X_{\Gamma_1}| >
n)$, it follows from Fact A that $\q^e(\Gamma_1>\tau_{n})$ decays
exponentially. For the second term of the right-hand side, beware
that
\begin{eqnarray*}
&&
   \q^e(\, \mbox{more than } \; n^2 \; \mbox{distinct sites are visited
   before} \; \tau_n)\\
&\le&
   \sum_{k=1}^{n} \q^e( \mbox{more than } \; n \;
   \mbox{distinct sites are visited at level}\, k)\,.
\end{eqnarray*}

\noindent If we denote by $t_i^k$ the first time when the $i$-th
distinct site of level $k$ is visited, we have, by the strong Markov
property,
\begin{eqnarray*}
       P_{\omega}^e\left( \mbox{more than } \; n \; \mbox{sites are visited at level}\; k\right)
&=&
       P_{\omega}^e\left(t_n^k<\infty\right)\\
&\le&
       P_{\omega}^e\left(t_{n-1}^k<\infty ,\, D\left(X_{t_{n-1}^k}\right)<\infty\right)\\
&=&
       E_{\omega}^e\left[\ind_{\{t_{n-1}^k<\infty\}}\left(1-\beta(X_{t_{n-1}^k})\right)\right]\,.
\end{eqnarray*}

\noindent The independence of the environments entails that
$$
E_{\q^e}\left[\ind_{\{t_{n-1}^k<\infty\}}\left(1-\beta(X_{t_{n-1}^k})\right)\right]=\q^e\left(t_{n-1}^k<\infty\right)E_{\Q}[1-\beta]\,.
$$

\noindent Consequently,
\begin{eqnarray}
    \q^e\left(t_n^k<\infty\right)
&\le&
    \q^e\left(t_{n-1}^k<\infty\right)E_{\Q}[1-\beta] \nonumber\\
&\le&
    \left(E_{\Q}[1-\beta]\right)^{n-1}\,, \label{LD:expsmall}
\end{eqnarray}

\noindent which leads to
\begin{eqnarray}
\q^e\left( \mbox{more than } \; n^2 \; \mbox{sites are visited
before}\; \tau_n\right) \le n\left(E_{\Q}[1-\beta]\right)^{n-1}\,,
\end{eqnarray}

\noindent which is exponentially small. We remark, for later use,
that equation (\ref{LD:expsmall}) holds without the assumption
$q_1=0$.
 For the last term of equation (\ref{LD:regen1}), we write
\begin{eqnarray*}
\q^e\left(\exists \, k\le n^2\,:\,N(\pi_k)> n \right) \le
\sum_{k=1}^{n^2}\q^e\left(N(\pi_k)> n\right)\,.
\end{eqnarray*}

\noindent Let $U:=\bigcup_{n\ge 0} (\bN^*)^n$ be the set of words,
where $(\bN)^0:=\{\emptyset\}$. Each vertex $x$ of $\t$ is naturally
associated with a word of $U$, and $\t$ is then a subset of $U$ (see
\cite{Ne} for a more complete description). For any $k\ge 1$,
\begin{eqnarray*}
      \q^e(N(\pi_k)> n)
&=&
      \sum_{x\in U}\q^e\left(x\in \t, \, N(x)>n,\,x=\pi_k\right)\\
&\le&
      \sum_{x\in U}E_{\Q}\left[\ind_{\{x\in \t\}}P_{\omega}^e(x=\pi_k)(1-\gamma(x))^{n}\right]\,,
\end{eqnarray*}

\noindent with the notation of (\ref{LD:gamma}). By independence,
\begin{eqnarray*}
                      \q^e(N(\pi_k)>n)
&\le&                 \sum_{x\in U}E_{\Q}\left[\ind_{\{x\in \t\}}P_{\omega}^e(x=\pi_k)\right]E_{\Q}\left[(1-\gamma(e))^{n}\right]\\
&=&                 E_{\Q}[(1-\gamma(e))^{n}]\,.
\end{eqnarray*}

\noindent Apply Lemma \ref{LD:beta} to complete the proof. $\Box$\\

\noindent {\it Proof of Lemma \ref{LD:beta}}. Let $\mu>0$ be such
that $q:=\Q(\beta(e)> \mu)>0$, and write
$$
R:=\inf \{ k \ge 1 \, : \, \exists |x|=k, \, \beta(x) \ge \mu \}\,.
$$

\noindent Let $x_R$ be such that $|x_R|=R$ and $\beta(x_R)\ge \mu$
and we suppose for simplicity that $x_R$ is a descendant of $e_1$.
We see that $\gamma(e)\ge \omega(e,e_1)\beta(e_1)\ge {c_1 \over
\nu(e)}\beta(e_1)$ by equation (\ref{LD:omega}). In turn, equation
(2.1) of \cite{aid07} implies that for any vertex $x$, we have
$$
{1\over \beta(x)} =1+{1\over \sum_{i=1}^{\nu(x)}
A(x_i)\beta(x_i)}\le 1 + {1\over \mbox{ess inf }A}{1\over
\beta(x_i)}\,,
$$

\noindent for any $1\le i\le \nu(x)$. By recurrence on the path from
$e_1$ to $x_R$, this leads to
$$
{1\over \beta(e_1)} \le 1 + {1\over \mbox{ess inf }A} +\ldots +
\left({1\over \mbox{ess inf }A}\right)^{R-1}{1\over \mu}\,.
$$

\noindent We deduce the existence of constants $c_4,\,c_5>0$ such
that
\begin{equation}
\label{LD:gammaest} {\gamma(e)}\ge {c_4\over \nu(e)}e^{-c_5R}\,.
\end{equation}

\noindent It yields that
\begin{eqnarray*}
        E_{\Q}\left[\left(1-\gamma(e)\right)^n\ind_{\{\nu(e)<\sqrt{n}\}}\right]
\le
        \Q\left(R> {1\over 4 c_5} \ln(n)\right) + e^{-n^{1/4+o(1)}}\,.
\end{eqnarray*}

\noindent We observe that
\begin{eqnarray*}
       \Q\left(R> {1 \over 4 c_5}
       \ln(n)\right)
\le
       \Q\left(\forall |x|= {1 \over 4c_5} \ln(n),\,\beta(x)>\mu\right)\,.
\end{eqnarray*}

\noindent By assumption, $q_1=0$; thus $\#\{x\in \t\,:\, |x|={1\over 4c_5} \ln(n)\}\ge 2^{1/4c_5\ln(n)} = :n^{c_6}$. As a consequence, $\Q\left(\forall |x|= {1\over 4c_5} \ln(n),\,\beta(x)>\mu\right)\le q^{n^{c_6}}$.
Hence, the proof of our lemma is reduced to find a stretched exponential bound for $E_{\Q}\left[\left(1-\gamma(e)\right)^n\ind_{\{\nu(e)\ge \sqrt{n}\}}\right]$.
For any $x\in \t$,
denote by $V_x^{\mu}$ the number of children $x_i$ of $x$ such that
$\beta(x_i)>\mu$. For $\varepsilon \in (0, \Q(\beta(e)>\mu))$,
\begin{eqnarray*}
&&
        E_{\Q}\left[ (1-\gamma(e))^n\ind_{\{\nu(e)\ge
        \sqrt{n}\}} \right]\\
&\le&
        \q^e\left( \nu(e)\ge \sqrt{n},\, V_{e}^{\mu} <
        \varepsilon \nu(e) \right) + E_{\Q}\left[ (1-\gamma(e))^n \ind_{\{
        V_{e}^{\mu}
        \ge \varepsilon \nu(e) \}}\right]\,.
\end{eqnarray*}

\noindent We apply Cram\'er's Theorem to handle with the first term on the right-hand side. Turning to the second one, the bound is clear once we observe the general inequality,
\begin{eqnarray}
     \gamma(e) &=&
      \sum_{k=1}^{\nu(e)}\omega(e,e_k)\beta(e_k)
\ge
      {c_1\over \nu(e)}\sum_{k=1}^{\nu(e)}\beta(e_k)
\ge
      {c_1 \mu \over \nu(e)}V_e^{\mu}\label{LD:gammasqrt}\,,
\end{eqnarray}

\noindent which is greater than $c_1\mu\varepsilon$ on $\{V_{e}^{\mu}
        \ge \varepsilon \nu(e)\}$. $\Box$\\

\noindent {\bf Remark 2.3}. As a by-product, we obtain that
$E_{\Q}\left[ (1-\gamma(e))^n\ind_{\{\nu(e)\ge
        \sqrt{n}\}} \right]\le e^{-n^{c_3}}$ without the assumption $q_1=0$.\\

\noindent{\it Proof of Proposition \ref{LD:timereg} : upper bound in
the case $s<1$.} We follow the strategy of the case $``q_1=0"$. The
proof boils down to the estimate of
\begin{eqnarray*}
&&\q^e(N(\pi_k)> n,\, D(e)=\infty) \label{LD:split}\\
&=&
\q^e(N(\pi_k)> n,\,\nu(\pi_k)<\sqrt{n},\,D(e)=\infty)+
\q^e(N(\pi_k)> n,\,\nu(\pi_k)\ge \sqrt{n},\,D(e)=\infty)\,.
\end{eqnarray*}

\noindent Let $x\in \t$ and consider the RWRE $(X_n,\,n\ge 0)$ when
starting from ${\buildrel \leftarrow \over x}$. Inspired by Lyons et
al. \cite{lpp96},  we propose to couple it with a random walk
$(Y_n'',\, n\ge 0)$ on $\bZ$. We first define $X_n''$ as the
restriction of $X_n$ on the path $[\![{\buildrel \leftarrow \over
e},x]\!]$. Beware that $X_n''$ exists only up to a time $T$, which
corresponds to the time when the walk $(X_n,\, n\ge 0)$ escapes the
path $[\![{\buildrel \leftarrow \over e},x]\!]$, id est leaves the
path and never comes back to it. After this time, we set $
X_n''=\Delta$ for some $\Delta$ in some space $\mathcal{E}$. Then $(
X_n '')_{n\ge 0}$ is a random walk on $[\![{\buildrel \leftarrow
\over e},x]\!]\cup \{\Delta\}$, whose transition probabilities are,
if $y\notin \{{\buildrel \leftarrow \over e},\, x,\,\Delta\}$,
\begin{eqnarray*}
    P_{\omega}^{{\buildrel \leftarrow \over x}}( X_{n+1}''=y_+\,|\,  X_n''=y)
&=&
    { \omega(y,y_+) \over \omega(y,y_+) + \omega(y,{\buildrel \leftarrow \over y})
    +\sum_{y_k\neq y_+}\omega(y,y_k)\beta(y_k)}\,,
    \\
    P_{\omega}^{{\buildrel \leftarrow \over x}}( X_{n+1}''={\buildrel \leftarrow \over y}\,|\,  X_n''=y)
&=&
    { \omega(y,{\buildrel \leftarrow \over y}) \over \omega(y,y_+) + \omega(y,{\buildrel \leftarrow \over y})
    +\sum_{y_k\neq y_+}\omega(y,y_k)\beta(y_k)}\,,
    \\
    P_{\omega}^{{\buildrel \leftarrow \over x}}(X_{n+1}''=\Delta\,|\,  X_n''=y)
&=&
    { \sum_{k=1}^{\nu(y)}\omega(y,y_k)\beta(y_k) \over \omega(y,y_+) + \omega(y,{\buildrel \leftarrow \over y})
    +\sum_{y_k\neq y_+}\omega(y,y_k)\beta(y_k)}\,,
\end{eqnarray*}

\noindent where $y_+$ is the child of $y$ which lies on the path
$[\![{\buildrel \leftarrow \over e},x]\!]$. Besides, the walk is
absorbed on $\Delta$ and reflected on ${\buildrel \leftarrow \over
e}$ and $x$. We recall that $s  :=   \mbox{ess sup} \,A$. We
construct the adequate coupling with a biased random walk
$(Y_n'')_{n\ge 0}$ on $\bZ$, starting from $|x|-1$, increasing with
probability $s/(1+s)$, decreasing otherwise and such that  $Y_n''
\ge |X_n''|$ as long as $X_n'' \neq \Delta$ (which is always
possible since $P_{\omega}( X_{n+1}''=y_+\,|\,  X_n''=y) \le {s\over
1+s}$). After time $T$, we let $Y_n$ move independently. By coupling
and then by gambler's ruin method, it leads to
\begin{eqnarray*}
     P_{\omega}^{{\buildrel \leftarrow \over x}}(T_x< T_{{\buildrel
     \leftarrow \over e}})
\le P_{\omega}^{|x|-1}(\,\exists \, n\ge 0\,:\,Y_n''=|x|\,)=s\,.
\end{eqnarray*}

\noindent It follows that
\begin{eqnarray*}
      1-P_{\omega}^x(T_x^*<T_{{\buildrel \leftarrow \over e}})
\ge
     \omega(x,{\buildrel \leftarrow \over x})\left(1-P_{\omega}^{{\buildrel \leftarrow \over x}}(T_x< T_{{\buildrel
     \leftarrow \over e}})\right)
\ge
     {c_{1}(1-s) \over \nu(x)}\,,
\end{eqnarray*}

\noindent by equation (\ref{LD:omega}). Hence,
\begin{eqnarray*}
&&       \q^e(N(\pi_k)>n,\,\nu(\pi_k) \le \sqrt{n},\,D(e)=\infty)\\
&=&
       \sum_{x\in U} E_{\Q}\left[ \ind_{\{\nu(x) \le \sqrt{n}\}}P_{\omega}^e\left(x=\pi_k, \, D(e)>T_x\right)P_{\omega}^x\left(N(x)>n, D(e)=\infty\right)
       \right]\\
&\le&
       \sum_{x\in U} E_{\Q}\left[P_{\omega}^e(x=\pi_k)\left( 1-{c_1(1-s)\over \sqrt{n}}\right)^n
       \right]=
       \left( 1-{c_1(1-s)\over \sqrt{n}}\right)^n\,,
\end{eqnarray*}

\noindent which decays stretched exponentially. On the other hand,
\begin{eqnarray*}
&&
        \q^e(N(\pi_k)>n,\,\nu(\pi_k) \ge \sqrt{n} , \, D(e)=\infty )\\
&\le&
        \q^e\left(\nu(\pi_k)\ge \sqrt{n} ,\, V_{\pi_k}^{\mu}< \varepsilon \nu(\pi_k)\right) +
        \q^e\left(N(\pi_k)>n ,\, V_{\pi_k}^{\mu} \ge \varepsilon \nu(\pi_k)\right)\,.
\end{eqnarray*}

\noindent with the notation introduced in the proof of Lemma
\ref{LD:beta}. We have
\begin{eqnarray*}
        \q^e\left(\nu(\pi_k)\ge \sqrt{n} ,\, V_{\pi_k}^{\mu}< \varepsilon \nu(\pi_k)\right)
&=&
        \Q\left(\nu(e)\ge \sqrt{n} ,\, V_{e}^{\mu}< \varepsilon \nu(e)\right)\,,
\end{eqnarray*}

\noindent which is stretched exponential by Cram\'er's Theorem. We
also observe that
\begin{eqnarray*}
        \q^e\left(N(\pi_k)>n ,\, V_{\pi_k}^{\mu} \ge \varepsilon \nu(\pi_k)\right)
&\le&
        E_{\q^e}\left[ \ind_{\{ V_{\pi_k}^{\mu} \ge \varepsilon \nu(x) \}} (1-\gamma(\pi_k))^n \right]\\
&=&
        E_{\Q}\left[ \ind_{\{ V_{e}^{\mu} \ge \varepsilon \nu(x) \}} (1-\gamma(e))^n \right]\le \left(1-c\mu\varepsilon\right)^n\,,
\end{eqnarray*}

\noindent by equation (\ref{LD:gammasqrt}). This completes the
proof. $\Box$

\subsection{The case $\Lambda<\infty$}

In this part, we suppose that $\Lambda<\infty$, where $\Lambda$ is defined by
$$
    \Lambda:=Leb \left\{t\in \r: \E[A^t]\le
    {1\over q_1}\right\}.
$$

\noindent We prove that the tail distribution of $\Gamma_1$ is
polynomial.
\begin{proposition}
\label{LD:timepol} If $\Lambda<\infty$, then
\begin{equation}
\lim_{n\rightarrow\infty}{1\over
\ln(n)}\ln\left(\S^e(\Gamma_1>n)\right)= -\Lambda\,.
\label{LD:timepoleq}
\end{equation}
\end{proposition}

\noindent {\it Proof of Proposition \ref{LD:timepol}}. Lemma 3.3 of
\cite{aid07} already gives
\begin{eqnarray*}
\liminf_{n\rightarrow\infty} { 1\over \ln(n)} \ln\left(\S^e(\Gamma_1>n)\right)\ge
-\Lambda.
\end{eqnarray*}

\noindent Hence, the lower bound of (\ref{LD:timepoleq}) is known.
The rest of the section is dedicated to the proof
of the upper bound.\\

\noindent We start with three preliminary lemmas. We first prove an
estimate for one-dimensional RWRE, that will be useful later on.
Denote by $(R_n,\,n\ge 0)$ a generic RWRE on $\bZ$ such that the
random variables $A(i)$, $i\ge 0$ are independent and have the
distribution of $A$, when we set for $i\ge 0$,
$$
A(i):={\omega_R(i,i+1)\over \omega_R(i,i-1)}
$$
with $\omega_R(y,z)$ the quenched probability to jump from $y$ to $z$. We denote by $P_{\omega,R}^k$ the quenched distribution associated with $(R_n,\,n\ge 0)$ when starting from $k$, and by $\P_R$ the distribution of the environment $\omega_R$. Let $c_7\in (0,1)$ be a constant whose value will be given later on. For any $k\ge \ell\ge 0$ and $n\ge 0$, we introduce the notation
\begin{eqnarray}
    p(\ell,k,n)
:=
    E_{\P_R}\left[ (1-c_7P_{\omega,R}^{\ell}(T_{\ell}^*
         >T_{0}\land T_{k}))^{n} \right]\,. \label{LD:defp}
\end{eqnarray}

\begin{lemma}
\label{LD:estone} Let $0<r<1$, and $\Lambda_r:=Leb \left\{t\in \r:
\E[A^t]\le
    {1\over r}\right\}$. Then, for any $\varepsilon>0$, we have for $n$ large enough,
$$
\sum_{k\ge \ell \ge 0}r^kp(\ell,k,n)\le n^{-\Lambda_r+\varepsilon}\,.
$$
\end{lemma}
\noindent {\it Proof}. The method used is very similar to that of
Lemma 5.1 in \cite{aid07}. We feel free to present a sketch of the
proof. We consider the one-dimensional RWRE $(R_n)_{n\ge 0}$. We
introduce for $k\ge \ell\ge 0$, the potential $V(0)=0$ and
\begin{eqnarray*}
  V(\ell)    &=& -\sum_{i=0}^{\ell-1}\ln(A(i))\,,\\
  H_1(\ell)  &=& \max_{0\le i\le \ell} V(i)-V(\ell)\,,\\
  H_2(\ell,k)&=& \max_{\ell\le i\le k} V(i)-V(\ell)\,.
\end{eqnarray*}

\noindent We know (e.g. \cite{zeitouni04}) that
\begin{eqnarray}
{e^{-H_2(\ell+1,k)}\over k+1}\le P_{\omega,R}^{\ell+1}\left(T_{k}<T_{\ell}\right)\le e^{-H_2(\ell+1,k)}\,,\\
{e^{-H_1(\ell)}\over k+1}\le
P_{\omega,R}^{\ell-1}\left(T_{-1}<T_{\ell}\right)\le e^{-H_1(\ell)}\,.
\end{eqnarray}

\noindent It yields that
\begin{eqnarray*}
P_{\omega,R}^{\ell}(T_{\ell}^*
>T_{0}\land T_{k}) \ge e^{-H_1(\ell)\land H_2(\ell,k)+ O(\ln k)}\,,
\end{eqnarray*}

\noindent where $O(\ln k)$ is a deterministic function. Let $\eta \in
(0,1)$.
\begin{eqnarray*}
p(\ell,k,n) &\le& (1-c_7n^{-1+\eta})^n + \P_R(H_1(\ell)\land
H_2(\ell,k)- O(\ln k) \ge (1-\eta)\ln(n))\\
&\le& e^{-c_8n^{\eta}} + \P_R(H_1(\ell)\land H_2(\ell,k)- O(\ln k) \ge
(1-\eta)\ln(n)) \,.
\end{eqnarray*}

\noindent In Section 8.1 of \cite{aid07}, we proved that for any
$s\in  (0,1)$, $E_{\P_R}\left[e^{\Lambda_s(H_1(\ell)\land
H_2(\ell,k))}\right] \le e^{k\ln(1/s)+o_s(k)}$, where $o_s(k)$ is
such that $o_s(k)/k $ tends to $0$ at infinity. This implies that,
defining $\widetilde{o}_s(k):=o_s(k)-\Lambda_s O(\ln k)$,
\begin{eqnarray*}
&&         s^k \P_R \left( H_1(\ell)\land H_2(\ell,k) -
         O(\ln k) \ge (1-\eta)\ln(n) \right)\\
&\le&
         s^k \left( 1 \land e^{k\ln(1/s) - \Lambda_s(1-\eta)\ln(n)+\widetilde{o}_s(k)}\right)\\
&\le&
         n^{-\Lambda_s(1-\eta)}\exp\left((k\ln(s)+\Lambda_s(1-\eta)\ln(n))\land \widetilde{o}_s(k)\right)
\,.
\end{eqnarray*}

\noindent Observe that there exists $M_s$ such that for any $k$ and
any $n$, we have $(k\ln(s) + \Lambda_s(1-\eta)\ln(n))\land
\widetilde{o}_s(k)\le \sup_{i\le M_s\ln(n)} \widetilde{o}(i) + \eta
\ln n$, and notice that $\sup_{i\le M_s\ln(n)} \widetilde{o}_s(i)$
is negligible towards $\ln(n)$. This leads to, for $n$ large enough,
$$
s^kp(\ell,k,n) \le s^ke^{-c_8n^{\eta}} + n^{-\Lambda_s(1-\eta)+2\eta}\,.
$$

\noindent Let $r\in (0,1)$ and $s>r$. We have
$$
r^kp(\ell,k,n) \le r^ke^{-c_8n^{\eta}} + \left(r \over s\right)^kn^{-\Lambda_s(1-\eta)+2\eta}\,.
$$

\noindent Lemma \ref{LD:estone} follows by choosing $\eta$ small enough and $s$ close enough to $r$. $\Box$\\

\noindent Let $Z_n$ represent the size of the $n$-th generation of the tree $\t$. We have the following result.
\begin{lemma}
\label{LD:betapol} There exists a constant $c_9>0$ such that for any
$H>0,\,B>0$ and $n$ large enough,
$$
E_{\Q}\left[\left(1-\gamma(e)\right)^{n}\ind_{\{Z_H>B\}}\right] \le
n^{-c_9B}\,.
$$
\end{lemma}
\noindent{\it Proof}. We have
\begin{eqnarray*}
E_{\Q}\left[\left(1-\gamma(e)\right)^{n}\ind_{\{Z_H>B\}}\right]
&\le& E_{\Q}\left[ (1-\gamma(e))^n\ind_{\{\nu(e)\ge
        \sqrt{n}\}} \right]
                           +
                           E_{\Q}\left[\left(1-\gamma(e)\right)^{n} \ind_{\{Z_H>B,\,\nu(e)\le \sqrt{n}\}}
                           \right]\\
        &\le&
e^{-n^{c_3}} + E_{\Q}\left[\left(1-\gamma(e)\right)^{n}
\ind_{\{Z_H>B,\,\nu(e)\le \sqrt{n}\}}
                           \right]
\end{eqnarray*}

\noindent by Remark 2.3. When $\nu(e)\le \sqrt{n} $, we have, by
(\ref{LD:gammaest}),
$$
{ \gamma(e)}\ge {c_4\over \sqrt{n}}e^{-c_5R}\,,
$$

\noindent with $R:=\inf \{ k \ge 1 \, : \, \exists |x|=k, \,
\beta(x) \ge \mu \}$ as before ($\mu>0$ is such that $q:=\Q(\beta(e)>\mu)>0$). Thus,
\begin{eqnarray*}
       E_{\Q}\left[\left(1-\gamma(e)\right)^{n} \ind_{\{Z_H>B,\,\nu(e)\le \sqrt{n}\}} \right]
\le
      \Q\left(R > {1\over 4c_5} \ln(n)+H,\, Z_H>B \right) + e^{-n^{1/4+o(1)}}\,.
\end{eqnarray*}

\noindent By considering the $Z_H$ subtrees rooted at each of the individuals in generation $H$, we see that
\begin{eqnarray*}
       \Q\left(R> c_{10} \ln(n)+H ,\, Z_H>B \right)
&=&
       E_{GW}\left[ \Q(R>c_{10}\ln(n))^{Z_H} \ind_{\{Z_H>B\}}\right]\\
&\le&
       \Q(R>c_{10}\ln(n))^{B}\,.
\end{eqnarray*}

\noindent If $R>c_{10}\ln(n)$, we have in particular $\beta(x)<\mu$ for each $|x|=c_{10}\ln(n)$ which implies that
\begin{eqnarray*}
       \Q\left(R> c_{10} \ln(n)+H ,\, Z_H>B \right)
\le
       E_{GW}\left[q^{Z_{c_{10} \ln(n)}}\right]^B\,.
\end{eqnarray*}

\noindent Let $t\in (q_1,1)$. For $n$ large enough, $E_{GW}\left[
q^{Z_{c_{10}\ln(n)}} \right] \le t^{c_{10}\ln(n)}=n^{{c_{10}\ln(t)}}
$, ($E_{GW}[q^{Z_n}]/q_1^n$ has a positive limit by Corollary 1 page
40 of \cite{AthNey72}). The lemma follows. $\Box$\\

Let $r\in (q_1,1)$, $\varepsilon>0$, $B$ be such that
\begin{equation}
\label{LD:defB} c_9B\varepsilon>2\Lambda
\end{equation}

\noindent and $H$
large enough so that
\begin{equation}
\label{LD:defH} GW( Z_H \le B) < r^{H}{1 \over B}<1.
\end{equation}

\noindent In particular, $c_{11}:=GW(Z_H>B)>0$.\\

Let $\nu(x,k)$ denote for any $x\in
\t$ the number of descendants of $x$ at generation $|x|+k$ ($\nu(x,1)=\nu(x)$), and let
\begin{eqnarray}
\mathcal{S}_H :=\{x\in \t\,:\, \nu(x,H)>B \}\,. \label{LD:tree1}
\end{eqnarray}

\noindent For any $x\in \t$, we call $F(x)$ the youngest ancestor of
$x$ which lies in $\mathcal{S}_H$, and $G(x)$ an oldest descendant
of $x$ in $\mathcal{S}_H$. For any $x,y\in \t$, we write $x\le y$ if
$y$ is a descendant of $x$ and $x<y$ if besides $x\neq y$. We define
for any $x \in \t$, $W(x)$ as the set of descendants $y$ of $x$ such
that there exists no vertex $z$ with $x< z \le y$ and $\nu(z,H)>B$.
In other words, $W(x)=\{ y \, : \, y \ge x,\,F(y)\le x \}$. We
define also
\begin{eqnarray*}
{\buildrel \circ \over W}(x)&:=& W(x)\backslash \{x\}\,, \\
{\partial W(x)} &:=& \{ y\,:\, {\buildrel \leftarrow \over y} \in
W(x),\, \nu(y,H)>B \}\,.
\end{eqnarray*}

\noindent Finally, let $W_j(e):=\{ x\,:\, |x|=j,\,x\in W(e) \}$.
\begin{lemma}
\label{LD:w} Recall that $m:=E_{GW}[\nu(e)]$ and $r$ is a real
belonging to $(q_1,1)$. We also recall that $H$ and $B$ verify
$GW(Z_H\le B)<r^H{1\over B}$. We have for any $j\ge 0$,
$$
    E_{GW}\left[W_{j}(e)\right]
< m\, r^{j-1}\,.
$$
\end{lemma}
\noindent {\it Proof}. We construct the subtree $\t_{H}$  of the
tree $\t$ by retaining only the generations $kH$, $k\ge 0$ of the
tree $\t$. Let
\begin{equation}
\mathbb{W}=\mathbb{W}(\t):=\{x\in\t_{H}: \forall y\in\t_{H}, (y<x)
\Rightarrow \nu(y,H)\le B \}\,. \label{LD:defw}
\end{equation}

\noindent The tree $\mathbb{W}$ is a Galton--Watson tree whose
offspring distribution is of mean $E_{GW}[Z_{H}\ind_{\{Z_{H}\le
B\}}]\le B\times GW( Z_H \le B) \le r^{H}$ by (\ref{LD:defH}). Then
for each child $e_i$ of $e$ (in the original tree $\t$), let
$\mathbb{W}^{i}:=\mathbb{W}(\t_{e_i})$ where $\t_{e_i}$ is the
subtree rooted at $e_i$. We conclude by observing that $W_j(e) \le
\sum_{i=1}^{\nu(e)}\#\{x\in \mathbb{W}^{i}\,:\,
|x|= 1+\lceil (j-1)/H \rceil\times H\}$ hence $E_{GW}\left[W_{j}(e)\right]\le E_{GW}[\nu(e)]r^{j-1}$. $\Box$\\

\noindent We still have $r\in (q_1,1)$ and $\varepsilon>0$. We prove
that for $n$ large enough, and $r$ and $\varepsilon$ close enough to
$q_1$ and $0$, we have
\begin{equation}
\label{LD:timepolbis} \q^e\left(\Gamma_1>n,\, D(e)=\infty\right)\le
c_{12}n^{-(1-2\varepsilon)\Lambda_r+ 3 \varepsilon}\,,
\end{equation}

\noindent where $\Lambda_r:=Leb\{t\in \r\,:\, \E[A^t]\le {1\over
r}\}$ as in Lemma \ref{LD:estone}. This suffices to prove
Proposition \ref{LD:timepol} since $\varepsilon$ and $\Lambda_r$ can
be arbitrarily close to $0$ and $\Lambda$, respectively. We recall
that we defined $B$, $H$ and $\mathcal{S}_H$ in
(\ref{LD:defB}),(\ref{LD:defH})
and (\ref{LD:tree1}).\\

The strategy is to divide the tree in subtrees in which vertices are constrained to have a small number of children (at most $B$ children at generation $H$). With $B=H=1$, we would have literally pipes. In general, the traps constructed are slightly larger than pipes. We then evaluate the time spent in such traps by comparison with a one-dimensional random walk.\\

\noindent We define $\pi_{k}^{\,s}$ as the $k$-th distinct site
visited in the set $\mathcal{S}_H$. We observe that
\begin{eqnarray}
&& \q^e\left(\Gamma_1>n,\,D(e)=\infty  \right) \label{LD:global}\\
\nonumber &\le& \q^e\left(\Gamma_1>\tau_{\ln^2(n)}\right) + \q^e\left(\mbox{more than
                   } \, \ln^4(n) \; \mbox{distinct sites are visited before}\; \tau_{\ln^2(n)}   \right)\\
&&   + \; \; \q^e\left(\exists \, k \le \ln^4(n),\,\exists \, x\in
W(\pi_k^{\,s}),\,
                  N(x)>n/\ln^4(n) \right) \nonumber \\
&& + \; \;
\q^e\left(\exists \, x\in
W(e),\,
                  N(x)>n/\ln^4(n),\,D(e)=\infty,\,Z_H\le B \right)\,. \nonumber
\end{eqnarray}

\noindent The first term on the right-hand side decays like
$e^{-\ln^2(n)}$ by Fact A, and so does the second term by equation
(\ref{LD:expsmall}). We proceed to estimate the third term on the
right-hand side of (\ref{LD:global}). Since
\begin{eqnarray*}
\q^e\left(\exists  k \le \ln^4(n),\exists  x\in
W(\pi_k^{\,s}),
                  N(x)>n/\ln^4(n) \right)
\le \sum_{k=1}^{\ln^4(n)}\q^e\left(\exists  x\in
W(\pi_k^{\,s}),
         N(x)>n/\ln^4(n)\right)
\end{eqnarray*}

\noindent we look at the rate of decay of $\q^e\left(\exists \, x\in
W(\pi_k^{\,s}),
         N(x)>n/\ln^4(n)\right)$ for any $k\ge 1$. We first show
         that the time spent at the frontier of $W(\pi_k^{\,s})$
         will be negligible. Precisely, we show
\begin{eqnarray}
\label{LD:c14} \q^e\left(N(\pi_k^{\,s})>n^{\varepsilon}\right)
        \le
c_{14} n^{-2\Lambda}\,,\\
\label{LD:c15}
        \q^e\left(\exists \, z\in \partial W(\pi_k^{\,s}),\, N(z)>n^{\varepsilon}\right)
\le
  c_{15}n^{-2\Lambda}\,.
\end{eqnarray}

\noindent As $ P_{\omega}^y( N(y)> n^{\varepsilon}) \le
      (1-\gamma(y))^{n^{\varepsilon}}$ for any $y\in \t$, we have,
\begin{eqnarray}
\nonumber \q^e\left(N(\pi_k^{\,s})>n^{\varepsilon}\right)
&=&
E_{\Q}\left[\sum_{\,y\in \, \mathcal{S}_H}P_{\omega}^e(\pi_k^{\,s}=y)P^{y}_{\omega}(N(y)>n^{\varepsilon})\right]\\
&\le&
E_{\Q}\left[\sum_{\,y\in\,\mathcal{S}_H}P_{\omega}^e(\pi_k^{\,s}=y)(1-\gamma(y))^{
n^{\varepsilon} } \label{LD:NF} \right]\,.
\end{eqnarray}

\noindent We would like to split the expectation
$E_{\Q}\left[P_{\omega}^e(\pi_k^{\,s}=y)(1-\gamma(y))^{
n^{\varepsilon} } \right]$ in two. However the random variable
$P_{\omega}^e(\pi_k^{\,s}=y)$ depends on the structure of the first
$H$ generations of the subtree rooted at $y$. Nevertheless, we are
going to show that, for some $c_{14}>0$,
$$
E_{\Q}\left[P_{\omega}^e(\pi_k^{\,s}=y)(1-\gamma(y))^{
n^{\varepsilon} } \right] \le
c_{14}E_{\Q}\left[P_{\omega}^e(\pi_k^{\,s}=y)\right]E_{\Q}\left[(1-\gamma(y))^{
n^{\varepsilon} }\,|\nu(y,H)>B \right]\,.
$$

\noindent Let $U:=\bigcup_{n\ge 0} (\bN^*)^n$ be, as before, the set
of words. We have seen that $U$ allows us to label the vertices of
any tree (see \cite{Ne}). Let $y\in U$ and let $\omega_y$ represent
the restriction of the environment $\omega$ to the outside of the
subtree rooted at $y$ (when $y$ belongs to the tree). For $1 \le L
\le H$, we denote by $y_L$ the ancestor of $y$ such that
$|y_L|=|y|-L$. We attach to each $y_L$ the variable
$\zeta(y_L):=\ind_{\{\nu(y_L,H)>B\}}$. We notice that there exists a
measurable function $f$ such that
$P_{\omega}^e(\pi_k^{\,s}=y)=f(\omega_y,\zeta)\ind_{\{\nu(y,H)>B\}}$
where $\zeta:=(\zeta(y_L))_{1\le L\le H}$. Let
 $\mathcal{E}(\omega_y):=\{e\in \{0,1\}^H\,:\,\Q(\zeta=e\,|\,\omega_y)>0\}$. We have
\begin{eqnarray*}
E_{\Q}\left[f(\omega_y,\zeta) \, | \, \omega_y\right]
\ge
\max_{e\,\in \,\mathcal{E}(\omega_y)} f(\omega_y,e)\Q\left(\zeta
\,=\,e\, | \, \omega_y \right)\,.
\end{eqnarray*}

\noindent We claim that there exists a constant $c_{13}>0$ such that
for almost every $\omega$ and any $e\in \mathcal{E}(\omega_y)$,
\begin{eqnarray*}
\Q\left(\zeta=e\, | \, \omega_y \right)\ge c_{13}\,.
\end{eqnarray*}

\noindent Let us prove the claim. If $\omega_y$ is such that $\nu({\buildrel \leftarrow
\over y})>B$, then $\mathcal{E}(\omega_y)=\{(1,\ldots,1)\}$ and
$\Q\left(\zeta=e\, | \, \omega_y \right)=1$. Therefore suppose
$\nu({\buildrel \leftarrow \over y}) \le B$ and let $h:= \max\{1\le
L \le H \,:\, \nu(y_L,L) \le B\}$. We observe that, for any $e\in
\mathcal{E}(\omega_y)$, we necessarily have $e_L=1$ for $h< L \le
H$. We are reduced to the study of
\begin{eqnarray*}
\Q\left(\zeta=e\,|\,\omega_y\right) = \Q\left(\bigcap_{1\,\le\,
L\,\le\,h}\{\zeta(y_L)=e_L\} \,\bigg|\, \omega_y \right)\,.
\end{eqnarray*}

\noindent For any tree $\mathcal{T}$, we denote by ${\mathcal{T}}^j$
the restriction to the $j$ first generations. Let also $\t_{y_h}$
designate the subtree rooted at $y_h$ in $\t$. Since
$\nu(y_{h},h)\le B$, we observe that $\t_{y_h}^h$ belongs almost
surely to a finite (deterministic) set in the space of all trees. We
construct the set
\begin{eqnarray*}
 \Psi(\t_{y_h}^h,e):=
 \{\mbox{tree}\;\mathcal{T} \,:\, \mathcal{T}^h=\t_{y_h}^h,\, GW(\mathcal{T}^{h+H})>0,\,\forall
 |x|\le
 2H,\, \nu_{\mathcal{T}}(x)\le B  \qquad \qquad && \\
  \,
  \forall\, 1\le L\le h,\, \nu_{\mathcal{T}}(y_L,h)>B  \; \mbox{if and only if}\;  e_L=1
 \}\,.&&
\end{eqnarray*}

\noindent We observe that $\Psi(\t_{y_K}^K,e)\neq \emptyset$ as soon
as $e\in \mathcal{E}(\omega_y)$. Let
$\widetilde{\Psi}(\t_{y_K}^K,e):=\{\mathcal{T}^{h+H},\,\mathcal{T}\in
\Psi(\t_{y_h}^h,e)\}$ be the same set but where the trees are restricted to the first $h+H$ generations. Since $\widetilde{\Psi}(\t_{y_K}^K,e)$ is
again included in a finite deterministic set in the space of trees,
we deduce that there exists $c_{13}>0$ such that, almost surely,
$$
\inf\{ GW(\mathcal{T}^{h+H}\,|\,\mathcal{T}^h),\,\mathcal{T} \in
\Psi(\t_{y_h}^h,e),\, e\in \mathcal{E}(\omega_y)\}\ge c_{13}\,.
$$

\noindent Consequently,
\begin{eqnarray*}
\Q\left(\zeta=e\, | \, \omega_y \right) \ge \Q(\t_{y_h}^{h+H} \in
\widetilde{\Psi}(\t_{y_h}^h,e)\,|\, \omega_y)\ge c_{13}\,,
\end{eqnarray*}

\noindent  as required. We get
\begin{eqnarray*}
E_{\Q}\left[f(\omega_y,\zeta) \, | \, \omega_y\right] \,\ge
\,c_{13}\,\max_{e\,\in \, \mathcal{E}(\omega_y)} f(\omega_y,e)
\,\ge
\,c_{13}\,f(\omega_y,\zeta)
\,.
\end{eqnarray*}

\noindent Finally we obtain, with $c_{14}:={1\over c_{13}}$,
$$
f(\omega_y,\zeta)\le
c_{14}\,E_{\Q}\left[f(\omega_y,\zeta) \, | \, \omega_y\right] \,.
$$

\noindent By (\ref{LD:NF}), it entails that
\begin{eqnarray*}
\q^e\left(N(\pi_k^{\,s})>n^{\varepsilon}\right)
&\le&
c_{14}\sum_{y\in U}E_{\Q}\left[ \ind_{\{\nu(y,H)>B\}} E_{\Q}\left[f(\omega_y,\zeta) \, | \, \omega_y\right]
(1-\gamma(y))^{ n^{\varepsilon} }\right]\\
&=&
c_{14}
\sum_{y\in U}E_{\Q}\left[f(\omega_y,\zeta)\right]
E_{\Q}\left[  \ind_{\{\nu(e,H)>B\}}(1-\gamma(e))^{ n^{\varepsilon} }\right]\\
&=&c_{14}
\sum_{y\in U}E_{\Q}\left[P_{\omega}^e(\pi_k^{\,s}=y)\right]
E_{\Q}\left[  (1-\gamma(e))^{ n^{\varepsilon} }\,|\, \nu(e,H)>B\right]\,.
\end{eqnarray*}

\noindent It implies that
\begin{eqnarray*}
\q^e\left(N(\pi_k^{\,s})>n^{\varepsilon}\right)
\le
     c_{14}E_{\Q}\left[(1-\gamma(e))^{n^{\varepsilon}}\, | \, Z_H>B\right] \le
c_{14} n^{-c_9\varepsilon B}\,,
\end{eqnarray*}

\noindent by Lemma \ref{LD:betapol}. Since $c_9\varepsilon
B>2\Lambda$, this leads to, for $n$ large,
$$
\q^e\left(N(\pi_k^{\,s})>n^{\varepsilon}\right)
        \le
c_{14} n^{-2\Lambda}
$$

\noindent which is equation (\ref{LD:c14}). Similarly, recalling
that $\partial W(y)$ designates the set of vertices $z$ such that
${\buildrel \leftarrow \over z}\in W(y)$ and $\nu(z,H)>B$, we have
that
\begin{eqnarray*}
   &&     \q^e\left(\exists \, y\in \partial W(\pi_k^{\,s}),\,
        N(y)>n^{\varepsilon}\right)\\
&\le&
        E_{\Q}\left[\sum_{y\in\mathcal{S}_H}P_{\omega}^e(\pi_k^{\,s}=y)\sum_{z\in
        \partial W(y)}(1-\gamma(z))^{n^{\varepsilon}}\right]\\
&\le&c_{14}
E_{\Q}\left[\sum_{y\in\mathcal{S}_H}P_{\omega}^e(\pi_k^{\,s}=y)\right]E_{GW}\left[\partial
W(e)\right]
E_{\Q}\left[(1-\gamma(e))^{n^{\varepsilon}}\,|\, Z_H>B\right]\\
&=&c_{14} E_{GW}\left[\partial
W(e)\right]E_{\Q}\left[(1-\gamma(e))^{n^{\varepsilon}}\,|\,
Z_H>B\right]\,.
\end{eqnarray*}

\noindent We notice that $E_{GW}[\partial W]\le
E_{GW}\left[\sum_{x\in W(e)} \nu(x)\right]=mE_{GW}\left[W(e)\right]$
which is finite by Lemma \ref{LD:w}. It yields, by Lemma
\ref{LD:betapol},
$$
        \q^e\left(\exists \, x\in W(\pi_k^{\,s}),\, N(G(x))>n^{\varepsilon}\right)
\le
  c_{15}n^{-2\Lambda}
$$

\noindent thus proving (\ref{LD:c15}). Our next step is then to find
an upper bound to the probability to spend most of our time at a
vertex $x$ belonging to some ${\buildrel \circ \over W}(y)$. To this
end, recall that $G(x)$ is an oldest descendant of $x$ such that
$\nu(x,H)>B$. We have just proved that the time spent at $y(=F(x))$
or $G(x)$ is negligible. Therefore, starting from $x$, the
probability to spend much time in $x$ is not far from the
probability to spend the same time without reaching $y$ neither
$G(x)$. Then, this probability is bound by coupling with a
one-dimensional random
walk.\\

\noindent Define $\widetilde{T}^{(\ell)}_x$ as the $\ell$-th time
the walk visits $x$ after visiting either $F(x)$ or $G(x)$, id est
$\widetilde{T}^{(1)}_x=T_x$ and,
$$
\widetilde{T}^{(\ell)}_x := \inf\{k> \widetilde{T}^{(\ell-1)}_x\,:\,
X_k=x,\,\exists \, i\in  (\widetilde{T}^{(\ell-1)}_x, k), \, X_i=F(x)
\;\mbox{or} \; G(x)\}\,.
$$

\noindent Let also
$N^{(\ell)}(x)=\sum_{k=\widetilde{T}^{(\ell)}(x)}^{\widetilde{T}^{(\ell+1)}(x)-1}\ind_{\{X_k=x\}}$
be the time spent at $x$ between $\widetilde{T}^{(\ell)}$ and
$\widetilde{T}^{(\ell+1)}$. We observe that, for any $k\ge 1$,
\begin{eqnarray}
&&
\q^e\left(\exists \, x\in W(\pi_k^{\,s}),
         N(x)>n/\ln^4(n)\right) \nonumber\\
&\le&
        \q^e\left(N(\pi_k^{\,s})>n^{\varepsilon}\right)
+
        \q^e\left(\exists \, x\in W(\pi_k^{\,s}),\,
        N(G(x))>n^{\varepsilon}\right) \nonumber \\
&&  \, +  \;         \q^e\left(\exists \, x\in {\buildrel \circ
\over W}(\pi_k^{\,s}),\,\exists \, \ell \le 2 n^{\varepsilon},\,
       N^{(\ell)}(x)>n^{1- 2 \varepsilon}\right) \nonumber \\
&\le& (c_{14}+c_{15})n^{-2\Lambda} + \sum_{\ell \le
2n^{\varepsilon}} \q^e\left(\exists \, x\in {\buildrel \circ \over
W}(\pi_k^{\,s}),\,
       N^{(\ell)}(x)>n^{1- 2 \varepsilon}\right)\,. \label{LD:glob2}
\end{eqnarray}

\noindent Since
\begin{eqnarray*}
\q^e( \exists \, x\in W(\pi_k^{\,s}),\,
       N^{(\ell)}(x)>n^{1- 2 \varepsilon})
\le
E_{\Q}\left[  \sum_{y\in\mathcal{S}_H}P_{\omega}^e(\pi_k^{\,s}=y)\sum_{x\in
        {\buildrel \circ \over W}(y)}P_{\omega}^x(
       N^{(\ell)}(x)>n^{1- 2 \varepsilon}) \right]\,,
\end{eqnarray*}

\noindent and by the strong Markov
property at $\widetilde{T}^{(\ell)}_x$,
\begin{eqnarray*}
        P_{\omega}^x\left(
        N^{(\ell)}(x)>n^{1- 2 \varepsilon}\right) \nonumber
&=&
        P_{\omega}^x\left({\widetilde T}^{(\ell)}_x<\infty\right)P_{\omega}^x\left(N^{(1)}(x)>n^{1 - 2
        \varepsilon}\right) \nonumber \\
&\le&
        P_{\omega}^x(N^{(1)}(x)>n^{1-2\varepsilon})\,,
\end{eqnarray*}

\noindent this yields
\begin{eqnarray}
&& \nonumber \q^e( \exists \, x\in W(\pi_k^{\,s}),\,
       N^{(\ell)}(x)>n^{1- 2 \varepsilon})\\
&\le&
 E_{\Q}\left[\sum_{y\in\mathcal{S}_H}P_{\omega}^e(\pi_k^{\,s}=y)\sum_{x\in
        {\buildrel \circ \over W}(y)}P_{\omega}^x(
       N^{(1)}(x)>n^{1- 2 \varepsilon})\right] \nonumber \\
&\le&
     c_{14}\,E_{\Q}\left[\sum_{y\in\mathcal{S}_H}P_{\omega}^e(\pi_k^{\,s}=y)\right]E_{\Q}\left[\sum_{x\in
        {\buildrel \circ \over W}(e)}P_{\omega}^x(N^{(1)}(x)>n^{1- 2 \varepsilon})\, \bigg | \,
        Z_H>B\right] \nonumber \\
&=&
        c_{14}\, E_{\Q}\left[\sum_{x\in
        {\buildrel \circ \over W}(e)}P_{\omega}^x(N^{(1)}(x)>n^{1- 2 \varepsilon})\, \bigg| \, Z_H>B\right]\,. \label{LD:eqfx}
\end{eqnarray}

\noindent  For any $x\in W(e)$, define, for any $y\in [\![e,G(x)]\!]$,
\begin{eqnarray*}
\widetilde{\omega}(y,y_+) &:=& {\omega(y,y_+) \over \omega(y,y_+) + \omega (y,{\buildrel \leftarrow \over y})}\,,\\
\widetilde{\omega}(y,{\buildrel \leftarrow \over y})
&:=& {\omega(y,{\buildrel \leftarrow \over y})\over \omega(y,y_+)+\omega(y,{\buildrel \leftarrow \over y})}\,,
\end{eqnarray*}

\noindent where as before $y_+$ represents the child of $y$ on the
path. We let $(\widetilde{X}_n)_{n\ge 0}$ be the random walk on
$[\![e,G(x)]\!]$ with the transition probabilities
$\widetilde{\omega}$ and we denote by
$\widetilde{P}_{\omega,x}(\cdot)$ the probability distribution of
$(\widetilde{X}_n,\,n\ge 0)$.  By Lemma 4.4 of \cite{aid07}, we have
the following comparisons:
\begin{eqnarray*}
        P^{{\buildrel \leftarrow \over x}}_{\omega}(T_x<T_{e})
&\le&
        \widetilde{P}^{{\buildrel \leftarrow \over x}}_{\omega,x}(T_x<T_{e})\,,\\
        P^{x_+}_{\omega}(T_{G(x)}<T_x)
&\le&
        \widetilde{P}^{x_+}_{\omega,x}(T_{G(x)}<T_{x})\,.
\end{eqnarray*}

\noindent Therefore,
\begin{eqnarray*}
&&
       P_{\omega}^x(T_x^*<T_{e}\land T_{G(x)})\\
&=&
       \omega(x,{\buildrel \leftarrow \over x})P_{\omega}^{{\buildrel \leftarrow \over
       x}}(T_x<T_{e})+\omega(x,x_+)P_{\omega}^{x_+}(T_x<T_{G(x)})+\sum_{i\le \nu(x):x_i\neq x^+}
\omega(x,x_i)(1-\beta(x_i))\\
&\le&
       \omega(x,{\buildrel \leftarrow \over x})\widetilde{P}_{\omega,x}^{{\buildrel \leftarrow \over x}}
(T_x<T_{e})+\omega(x,x_+)\widetilde{P}_{\omega,x}^{x_+}(T_x<T_{G(x)})
       + \sum_{i\le \nu(x):x_i\neq x_+}\omega(x,x_i)\\
&=&
       1-\left(\omega(x,{\buildrel \leftarrow \over x}) + \omega(x,x_+)
       \right)\widetilde{P}_{\omega,x}^{x}(T_x^* > T_{e}\land T_{G(x)})\,.
\end{eqnarray*}

\noindent Since $\nu(x) \le B$ (for $x \in {\buildrel \circ \over
W}(e)$), we find by (\ref{LD:omega}) a constant $c_{16} \in (0,1)$
such that $\omega(x,{\buildrel \leftarrow \over x}) +
\omega(x,x_+)\ge c_{16}$. It yields that
\begin{eqnarray*}
P_{\omega}^x(T_x^*<T_{e}\land T_{G(x)}) \le
1-c_{16}\widetilde{P}_{\omega,x}^x(T_x^*
>T_{e}\land T_{G(x)})\,.
\end{eqnarray*}

\noindent We observe that, for any $x\in W(e)$, with the notation of
(\ref{LD:defp}) and taking $c_7:=c_{16}$,
\begin{eqnarray*}
E_{\P}\left[\left(1-c_{16}\widetilde{P}_{\omega,x}^x(T_x^*
>T_{e}\land T_{G(x)})\right)^n\right]=p(|x|,|G(x)|,n)\,.
\end{eqnarray*}

\noindent It follows that
\begin{eqnarray*}
      E_{GW}\left[  \sum_{x\in
        {\buildrel \circ \over W}(e)}\p^x(N^{(1)}(x)>n^{1- 2 \varepsilon})
       \right]
\le
      E_{GW}\left[  \sum_{x\in
        {\buildrel \circ \over W}(e)}p(|x|,|G(x)|,n^{1-2\varepsilon})\right]\,.
\end{eqnarray*}

\noindent On the other hand, $\sum_{x\in
        W(e)}p(|x|,|G(x)|,n^{1-2\varepsilon}) \le \sum_{y\in \partial W(e)} \sum_{x\le y} p(|x|,|y|,n^{1-2\varepsilon})$.
It implies that
\begin{eqnarray*}
E_{GW}\left[  \sum_{x\in
        {\buildrel \circ \over W}(e)}\p^x(N^{(1)}(x)>n^{1- 2 \varepsilon})
       \right]
&\le&
      \sum_{j\ge 0}E_{GW}\left[\#\{y\in \partial W(e), |y|=j\}\right]\left(\sum_{i \le j}
p(i, j,
       n^{1-2\varepsilon})\right)\\
&\le&
       m\sum_{j\ge 0}E_{GW}\left[W_{j-1}(e)\right]\left(\sum_{i \le j} p(i, j,
       n^{1-2\varepsilon})\right)\,.
\end{eqnarray*}

\noindent By Lemmas \ref{LD:estone} and \ref{LD:w}, for $n$ large
enough,
\begin{eqnarray}
&&E_{GW}\left[  \sum_{x\in
        {\buildrel \circ \over W}(e)}\p^x(N^{(1)}(x)>n^{1- 2 \varepsilon})
       \right]
\le
       m^2\, \sum_{ j \ge 0 } r^{j-2} \left(\sum_{i \le j} p(i, j,
       n^{1-2\varepsilon})\right) \le  n^{-(1-2\varepsilon)\Lambda_r+\varepsilon}\,.\nonumber\\
\label{LD:c17}
\end{eqnarray}

\noindent  Supposing $r$ and $\varepsilon$ close enough to $q_1$ and
$0$, equation (\ref{LD:c17}) combined with (\ref{LD:glob2}) and
(\ref{LD:eqfx}), shows that, for any $k\ge 1$,
$$
\q^e\left(\exists \, x\in W(\pi_k^{\,s}),
         N(x)>n/\ln^4(n)\right)
\le c_{17}n^{-(1-2\varepsilon)\Lambda_r+2\varepsilon}\,.
$$

\noindent We arrive at
\begin{eqnarray}
\label{LD:part3} \q^e\left(\exists \, k \le \ln^4(n),\,\exists \,
x\in W(\pi_k^{\,s}),\,
                  N(x)>n/\ln^4(n) \right) \le c_{18}n^{-(1-2\varepsilon)\Lambda_r+3\varepsilon}\,.
\end{eqnarray}

\noindent Finally, the estimate of
$\q^e\left(\exists \, x\in
W(e),\,
                  N(x)>n/\ln^4(n),\,D(e)=\infty,\, Z_H\le B \right)
$ in (\ref{LD:global}) is similar. Indeed,
\begin{eqnarray*}
&&\q^e\left(\exists \, x\in
W(e),\,
                  N(x)>n/\ln^4(n),\,D(e)=\infty,\, Z_H\le B \right)\\
&\le&
        \q^e\left(N(e)>n^{\varepsilon},\,D(e)=\infty,\,\nu(e)\le B\right)
+
        \q^e\left(\exists \, x\in W(e),\,
        N(G(x))>n^{\varepsilon}\right) \nonumber\\
&&  \, +  \;         \q^e\left(\exists \, x\in W(e),\,\exists
\, \ell \le 2 n^{\varepsilon},\,
       N^{(\ell)}(x)>n^{1- 2 \varepsilon}\right)\,.
\end{eqnarray*}

\noindent We have
\begin{eqnarray*}
\q^e\left(N(e)>n^{\varepsilon},\,D(e)=\infty,\,\nu(e)\le B\right)
&\le&
E_{\Q}\left[(1-\omega(e,{\buildrel \leftarrow \over e}))^{n^{\varepsilon}}\ind_{\{\nu(e)\le B\}}\right]\\
&\le&
\left(1-c_{1}/B\right)^{n^{\varepsilon}}\,,
\end{eqnarray*}

\noindent by (\ref{LD:omega}). By equation (\ref{LD:c15}),
\begin{eqnarray*}
\q^e\left(\exists \, x\in W(\pi_k^{\,s}),\, N(G(x))>n^{\varepsilon}\right)
\le
  c_{15}n^{-2\Lambda}\,.
\end{eqnarray*}

\noindent Finally,
\begin{eqnarray*}
\q^e\left(\exists \, x\in {\buildrel \circ \over W}(e),\,\exists \,
\ell \le 2 n^{\varepsilon},\,
       N^{(\ell)}(x)>n^{1- 2 \varepsilon}\right)
&\le&          \sum_{\ell \le 2n^{\varepsilon}} \q^e\left(\exists \,
x\in {\buildrel \circ \over W}(e),\,
       N^{(\ell)}(x)>n^{1- 2 \varepsilon}\right)  \\
&\le&   2n^{\varepsilon} \q^e\left(\exists \, x\in {\buildrel \circ
\over W}(e),\,
       N^{(1)}(x)>n^{1- 2 \varepsilon}\right)\\
&\le&
      2n^{\varepsilon}E_{GW}\left[  \sum_{x\in
        {\buildrel \circ \over W}(e)}\p^x(N^{(1)}(x)>n^{1- 2 \varepsilon})
       \right]\\
&\le& c_{17}n^{-(1-2\varepsilon)\Lambda_r+2\varepsilon}\,,
\end{eqnarray*}

\noindent by (\ref{LD:c17}). We deduce that, for $n$ large enough,
\begin{eqnarray}
\label{LD:part4} \q^e\left(\exists \, x\in W(e),\,
                  N(x)>n/\ln^4(n),\,D(e)=\infty,\, Z_H\le B \right)
\le
n^{-(1-2\varepsilon)\Lambda_r+3\varepsilon}\,.
\end{eqnarray}

\noindent In view of $(\ref{LD:global})$ combined with
$(\ref{LD:part3})$ and $(\ref{LD:part4})$, equation
(\ref{LD:timepolbis}) is proved, and Proposition \ref{LD:timepol}
follows. $\Box$


\section{Large deviations principles}

We recall the definition of the first regeneration time
$$
\Gamma_1:=\inf\left\{k>0 \, : \nu(X_k)\ge 2,\,
D(X_k)=\infty,\,k=\tau_{|X_k|} \right\}.
$$

\noindent We define by iteration $$\Gamma_n:=\inf\left\{k>
\Gamma_{n-1} \, : \nu(X_k)\ge 2,\, D(X_k)=\infty
,\,k=\tau_{|X_k|}\right\}$$ for any $n\ge 2$. We have the following
fact (points (i) to (iii) are already discussed in \cite{aid07};
point (iv) is shown in \cite{gross} in the case of regular trees and
in \cite{lpp96} in the case of biased random walks, and is easily
adaptable to our case).

\bigskip

\noindent {\bf Fact B} {\it

{\rm (i)} For any $n\ge 1$,\, $\Gamma_n<\infty$ ~~$\q^e$-a.s.

{\rm (ii)} Under
$\q^e$,\,$(\Gamma_{n+1}-\Gamma_n,|X_{\Gamma_{n+1}}|-|X_{\Gamma_n}|),\,n\ge
1$ are independent and distributed as $(\Gamma_1,|X_{\Gamma_1}|)$
under the distribution $\S^e$.

{\rm (iii)} We have $E_{\S^e}[\,|X_{\Gamma_1}|\,]<\infty$.

{\rm (iv)} The speed $v$ verifies $v={E_{\S^e}[\,|X_{\Gamma_1}|\,]
\over E_{\S^e}[\,\Gamma_1\,]}$.
   }\\

\noindent The rest of the section is devoted to the proof of
Theorems \ref{LD:speedup} and \ref{LD:slowdown}. It is in fact
easier to prove them when conditioning on never returning to the
root. Our theorems become

\begin{theorem}
\label{LD:speedupcond} {\bf (Speed-up case)} There exist two
continuous, convex and strictly decreasing functions $I_a \le I_q$
from $[1,1/v]$ to $\mathbb{R}_+$, such that $I_a(1/v)=I_q(1/v)=0$
and for $ a<b$, $b\in [1, 1/v]$,
\begin{eqnarray}
     \lim_{n \rightarrow \infty} { 1 \over n }\ln\left(\q^e\left({\tau_n\over n}\in ]a, b]\,\bigg|\, D(e)=\infty\right)\right)
&=&
     - I_a(b)\,, \label{LD:devupaup2}\\
     \lim_{n\rightarrow\infty} {1\over n}\ln \left(P_{\omega}^e\left( {\tau_n \over n} \in ]a, b] \,\bigg|\, D(e)=\infty\right)\right)
&=&
     - I_q(b)\,. \label{LD:devupqup2}
\end{eqnarray}
\end{theorem}

\begin{theorem}
\label{LD:slowdowncond} {\bf (Slowdown case)} There exist two
continuous, convex functions $I_a \le I_q$ from $[1/v,+\infty[$ to
$\mathbb{R}_+$, such that $I_a(1/v)=I_q(1/v)=0$ and for any $1/v \le
a < b$,
\begin{eqnarray}
     \lim_{n\rightarrow\infty} {1\over n}\ln\left(\q^e\left({\tau_n \over n} \in [a,b[\,\bigg|\, D(e)=\infty\right)\right)
&=&
     -I_a(a)\,, \label{LD:devupaslow2}\\
     \lim_{n\rightarrow\infty} {1\over n}\ln \left(P_{\omega}^e\left({\tau_n \over n}\in [a,b[\,\bigg|\, D(e)=\infty\right) \right)
&=&
     -I_q(a)\,. \label{LD:devupqslow2}
\end{eqnarray}
If $\mbox{ess inf }A =:i> \nu_{min}^{-1}$, then $I_a$ and $I_q$ are
strictly increasing on $[1/v,+\infty[$. If $i\le \nu_{min}^{-1}$,
then $I_a=I_q=0$.
\end{theorem}

\bigskip
\bigskip

\noindent Theorems \ref{LD:speedup} and \ref{LD:slowdown} follow
from Theorems \ref{LD:speedupcond} and \ref{LD:slowdowncond} and the
following proposition.
\begin{proposition}
\label{LD:rootdev} We have, for $a<b\le 1/v$,
\begin{eqnarray}
\lim_{n\rightarrow\infty} {1\over n}\ln\left(\q^e({\tau_n\over n}
\in ]a,b]) \right) &=& \lim_{n\rightarrow\infty} {1\over
n}\ln\left(\q^e({\tau_n\over n}
\in ]a,b]\,|\, D(e)=\infty)\right)\,, \label{LD:unconspeedupa} \\
\lim_{n\rightarrow\infty} {1\over
n}\ln\left({P_{\omega}^e({\tau_n\over n} \in ]a,b])}\right)&=&
\lim_{n\rightarrow\infty} {1\over
n}\ln\left(P_{\omega}^e({\tau_n\over n} \in ]a,b] \,|\,
D(e)=\infty)\right) \,.\label{LD:unconspeedupq}
\end{eqnarray}

\noindent Similarly, in the slowdown case, we have for $1/v \le
a<b$,
\begin{eqnarray}
\lim_{n\rightarrow\infty} {1\over n}\ln\left({\q^e({\tau_n\over n} \in [a,b[)
}\right)&=&
\lim_{n\rightarrow\infty} {1\over n}\ln\left(\q^e({\tau_n \over n} \in [a,b[ \,|\, D(e)=\infty)\right)\,, \label{LD:unconslowdowna} \\
\lim_{n\rightarrow\infty} {1\over n}\ln\left({P_{\omega}^e({\tau_n
\over n} \in [a,b[)}\right) &=& \lim_{n\rightarrow\infty} {1\over
n}\ln\left(P_{\omega}^e({\tau_n \over n} \in [a,b[ \,|\,
D(e)=\infty) \right)\label{LD:unconslowdownq}\,.
\end{eqnarray}
\end{proposition}

Theorems \ref{LD:speedupcond} and \ref{LD:slowdowncond} are proved
in two distinct parts for sake of clarity. Proposition
\ref{LD:rootdev} is proved in subsection \ref{LD:sectrootdev}.

\subsection{Proof of Theorem \ref{LD:speedupcond}}
\label{LD:sectspeedupcond} For any real numbers $h\ge 0$ and $b \ge
1$, any integer $n\in \bN$ and any vertex $x\in \t$ with $|x|=n$,
define
\begin{eqnarray*}
A(h,b,x)&:=&\{\omega\,:\,P_{\omega}^e\left(\tau_n=T_x,\,\tau_n\le
bn,\,T_{{\buildrel \leftarrow \over e}}> \tau_n \right)\ge e^{-hn}\}\,,\\
e_n(h,b)&:=& E_{\Q}\left[\sum_{|x|=n}\ind_{A(h,b,x)}\right]\,.
\end{eqnarray*}

\noindent We define also for any $b\ge 1$
\begin{eqnarray*}
h_c(b)&:=& \inf\{h\ge 0\,:\, \exists\, p\in \bN, \, e_p(h,b)>0\}\,.
\end{eqnarray*}

\begin{lemma}
\label{LD:defe} There exists for any $b \ge 1$ and $h > h_c(b)$, a
real $e(h,b)> 0$ such that
$$
\lim_{n\rightarrow\infty}{1\over
n}\ln(e_n(h,b))=\ln(e(h,b))\,.
$$
Moreover, the function $(h,b)\rightarrow \ln(e(h ,b))$ from
$\{(h,b)\in \r_+\times [1,+\infty[\,:\, h>h_c(b)\}$ to $\r$ is
concave, is nondecreasing in $h$ and in $b$, and
$$
\lim_{h\rightarrow\infty} {\ln(e(h,b))}=\ln(m)\,.
$$
\end{lemma}
\noindent{\it Proof.} Let $x\le y$ be two vertices of $\t$ with
$|x|=n$ and $|y|=n+m$. We observe that
\begin{eqnarray*}
           A(h,b,y)
&\supset&
           A(h,b,x) \cap \{\omega\,:\,P_{\omega}^x(\tau_{n+m}=T_y,\,
           \tau_{n+m}\le bm,\,T_{{\buildrel \leftarrow \over x}}>\tau_{n+m})\ge
           e^{-hm}\}\\
&=:&
           A(h,b,x) \cap A_x(h,b,y).
\end{eqnarray*}

\noindent It yields that
\begin{eqnarray}
       e_{n+m}(h,b)
&\ge&
       E_{\Q}\left[\sum_{|x|=n}\ind_{A(h,b,x)}\sum_{|y|=n+m,y\ge
       x}\ind_{A_x(h,b,y)}\right]\nonumber \\
&=&
       E_{\Q}\left[\sum_{|x|=n}\ind_{A(h,b,x)}\right]E_{\Q}\left[\sum_{|x|=m}\ind_{A(h,b,x)}\right] \nonumber \\
&=&
      e_{n}(h,b)e_m(h,b)\label{LD:subad}\,.
\end{eqnarray}

\noindent Let $h>h_c$ and $p$ be such that $e_p(h_c,b)>0$, where we
write $h_c$ for $h_c(b)$. Then $e_{np}(h_c,b)>0$ for any $n\ge 1$.
We want to show that $e_{k}(h,b)>0$ for $k$ large enough. By
(\ref{LD:omega}), $\omega(e,e_1)\ge c_1$ if $\nu(e)=1$ so that
$e_k(-\ln(c_1),b)\ge q_1^k$. Let $n_c$ be such that
$e^{-h_cn_c}c_1\ge e^{-hn_c}$. We check as before that for any $n\ge
n_c$, and any $r\le p$, we have indeed
\begin{eqnarray*}
e_{np+r}(h,b) &\ge& e_{np}(h_c,b)e_{r}(-\ln(c_1),b)\\
&\ge& e_{np}(h_c,b)q_1^{r}>0\,.
\end{eqnarray*}

\noindent Thus (\ref{LD:subad}) implies that
\begin{eqnarray}
\lim_{n\rightarrow\infty}{1\over n}\ln(e_n(h,b))=\sup\left\{{1\over
k}\ln(e_k(h,b)),\,k\ge 1\right\}=:\ln(e(h,b))\,,\label{LD:supehb}
\end{eqnarray}

\noindent with $e(h,b)>0$. Similarly, we can check that
\begin{eqnarray*}
e_n(t h_1 + (1-t)h_2,t b_1 +(1-t)b_2)\ge
e_{nt}(h_1,b_1)e_{n(1-t)}(h_2,b_2)\,,
\end{eqnarray*}

\noindent which leads to
\begin{eqnarray*}
\ln(e(t h_1 + (1-t)h_2,t b_1 +(1-t)b_2))\ge
t\ln(e(h_1,b_1))+(1-t)\ln(e(h_2,b_2))\,,
\end{eqnarray*}

\noindent hence the concavity of $(h,b)\rightarrow \ln(e(h,b))$. The
fact that $e(h,b)$  is nondecreasing in $h$ and in $b$ is direct.
Finally, $\limsup_{h\rightarrow \infty}\ln(e(h,b))\le \ln(m)$ and $
\liminf_{h\rightarrow \infty}{\ln(e(h,b))}\ge \liminf_{h\rightarrow
\infty}{\ln(e_1(h,b))}=\ln(m)$ by dominated convergence. $\Box$\\

In the rest of the section, we extend $e(h,b)$ to
$\r_+\times[1,+\infty[$ by taking $e(h,b)=0$ for $h\le h_c(b)$.\\

\begin{corollary}
\label{LD:s} Let $S:=\{h\ge 0: e(h,b)>1\}$ and $S':=\{h\ge 0:
e(h,b)\ge 1\}$. We have
$$
\sup\{e^{-h}\,e(h,b),\,h\in S\}=\sup\{e^{-h}\,e(h,b),\,h\in S'\}\,.
$$
\end{corollary}
\noindent{\it Proof.} Let $M:=\inf\{h\,:\,e(h,b)>1\}$. We claim that
if $h<M$, then $e(h,b)<1$. Indeed, suppose that there exists $h_0<M$
such that $e(h_0,b)\ge 1$. Then $e(h_0,b)=1$ by definition of $M$,
so that $e(h,b)$ is constant equal to $1$ on $[h_0,M[$. By
concavity, $\ln(e(h,b))$ is equal to $0$ on $[h_0,+\infty[$, which
is impossible since it tends to $\ln(m)$ at infinity.
The corollary follows. $\Box$\\

We have the tools to prove Theorem \ref{LD:speedup}.\\ \\
{\it Proof of Theorem \ref{LD:speedup}.} For $b \in [1,+\infty[$,
let
\begin{eqnarray*}
     J_a(b)
&:=&
     -\sup\{-h + \ln(e(h,b))\,,\,h\ge 0\}\,,\\
     J_q(b)
&:=&
     -\sup\{-h + \ln(e(h,b))\,,\,h\in S\}\,.
\end{eqnarray*}

\noindent Define then for any $b\le1/v$,
\begin{eqnarray*}
I_a(b) &=& J_a(b)\,,\\
I_q(b) &=& J_q(b)\,.
\end{eqnarray*}

\noindent We immediately see that $I_a \le I_q$. The convexity of
$J_a$ and $J_q$ stems from the convexity of the function
$h-\ln(e(h,b))$. Indeed, let $J$ represent either $J_a$ or $J_q$ and
let $1 \le b_1 \le b_2$ and $t\in [0,1]$. Denote by $h_1$, $h_2$,
$b$ and $h$ the reals that verify
\begin{eqnarray*}
J(b_1) &=& h_1-\ln(e(h_1,b_1))\,,\\
J(b_2) &=& h_2-\ln(e(h_2,b_2))\,,\\
h &:=& th_1+(1-t)h_2\,,\\
b &:=& tb_1 + (1-t)b_2\,.
\end{eqnarray*}

\noindent We observe that
\begin{eqnarray*}
J(b) &\le& h-\ln(e(h,b))\\
 &\le&
t(h_1-\ln(e(h_1,b_1)))+(1-t)(h_2-\ln(e(h_2,b_2)))=tJ(b_1)+(1-t)J(b_2)\,
\end{eqnarray*}

\noindent which proves the convexity. We show now that, for any
$b\ge 1$,
\begin{eqnarray}
\lim_{n\rightarrow\infty}{1\over n}\ln\left( \q^e\left(\tau_n<T_{{\buildrel \leftarrow \over e}}, \,  \tau_n\le b n\right)  \right)
&=& -J_a(b)\,, \label{LD:ia} \\
\lim_{n\rightarrow\infty}{1\over n}\ln\left(
P_{\omega}^e\left(\tau_n<T_{{\buildrel \leftarrow \over e}}, \,
\tau_n\le b n\right)  \right) &=& -J_q(b)  \,.  \label{LD:iq}
\end{eqnarray}

\noindent We first prove (\ref{LD:ia}). Since
$\q^e\left(\tau_n<T_{{\buildrel \leftarrow \over e}}, \,  \tau_n\le
b n\right) \ge
      e^{-hn}e_n(h,b)$ for any $h\ge 0$, we have

\begin{eqnarray*}
     \liminf_{n\rightarrow\infty} {1\over n}\ln \left( \q^e(\tau_n < T_{{\buildrel
     \leftarrow \over e}},  \, \tau_n\le bn)\right)
\ge
     -I_a(b).
\end{eqnarray*}

\noindent Turning to the upper bound, take a positive integer $k$.
We observe that
\begin{eqnarray*}
       \q^e\left(\tau_n<T_{{\buildrel \leftarrow \over e}}, \, \tau_{n}\le
       bn\right)
&\le&
       \sum_{\ell=0}^{k-1}
       e^{- n \ell /k}e_{n}\left((\ell +1)/k,b\right)\\
&\le&
       ke^{n/k}\sup\{e^{-hn}e_n(h,b) , \, h\ge 0 \}\,.
\end{eqnarray*}

\noindent Therefore,
\begin{eqnarray*}
           \limsup_{n\rightarrow\infty} {1\over n}\ln\left(\q^e\left(\tau_n< T_{{\buildrel
           \leftarrow \over e}} , \,\tau_n\le bn \right) \right)
\le
           {1\over k} - J_a(b)\,.
\end{eqnarray*}

\noindent Letting $k$ tend to infinity gives the upper bound of
(\ref{LD:ia}).\\

To prove equation (\ref{LD:iq}), let $k$ be still a positive integer
and $ h \in S$. Denote by $V_{pk}(\t)$ the set of vertices $ |x| =
pk $ such that $ P_{\omega}^{x_{\ell-1}}\left( \tau_{ \ell
k}<T_{{\buildrel \leftarrow \over x}_{\ell-1}}, \,\tau_{ \ell
k}=T_{x_{\ell}} \le b k \right) \ge e^{-hk} $ for any $\ell \le p$,
where  $x_{\ell}$ represents the ancestor of $x$ at generation $\ell
k$. Call $ V ( \t ):=\cup_{p\ge 0}V_{pk}(\t) $ the subtree thus
obtained. We observe that $ V $ is a Galton--Watson tree of mean
offspring $ e_k(h,b) $. Let
\begin{eqnarray*}
    \mathcal{T}_{k,h}
:=
    \{ \t\,:\, V(\t)\,\mbox{is infinite} \} \, .
\end{eqnarray*}

\noindent Take $\t \in \mathcal{T}_{k,h}$. For any $x\in V_{pk}$, we
have
\begin{eqnarray*}
 &&       P_{\omega}^e\left( \tau_{p k}<T_{{\buildrel \leftarrow \over e}},
        \,\tau_{p k}=T_x \le b p k \right)\\
 &\ge&
        P_{\omega}^e\left( \tau_{k}<T_{{\buildrel \leftarrow \over e}},
        \,\tau_{k}=T_{x_1} \le bk \right)
        \ldots
        P_{\omega}^{x_{k-1}}\left( \tau_{p k}<T_{{\buildrel \leftarrow \over x_{k-1}}},
        \,\tau_{p k}=T_x \le bk \right)
 \ge
        e^{-hpk}\,.
\end{eqnarray*}

\noindent It implies that
\begin{eqnarray*}
      P_{\omega}^e\left( \tau_{p k}<T_{{\buildrel \leftarrow \over e}},
      \,\tau_{p k} \le b p k \right)
\ge
      e^{-h p k}\#V_{p k}(\t)\,.
\end{eqnarray*}

\noindent By the Seneta--Heyde Theorem (see \cite{AthNey72} page 30
Theorem 3),
$$
\lim_{p\rightarrow\infty}{1\over p}{\ln\left( \#V_{p k}(\t) \right)}
= \ln(e_k(h,b)) \qquad \Q-\mbox{a.s.}
$$

\noindent It follows that, as long as $\t \in \mathcal{T}_{k,h}$,
\begin{eqnarray*}
     \liminf_{p\rightarrow\infty} {1\over p k}\ln\left( P_{\omega}^e\left( \tau_{p k}<T_{{\buildrel \leftarrow \over e}},
      \,\tau_{p k} \le b p k \right)
     \right)
\ge
     -h + {1\over k}\ln(e_k(h,b))\,.
\end{eqnarray*}

\noindent Notice that
\begin{eqnarray*}
      P_{\omega}^e\left( \tau_{n}<T_{{\buildrel \leftarrow \over e}},
      \,\tau_{n} \le bn \right)
\ge
     P_{\omega}^e\left( \tau_{p k}<T_{{\buildrel \leftarrow \over e}},
      \,\tau_{p k} \le bp k \right)\min_{|x|=p k}P_{\omega}^x\left( \tau_{n}<T_{{\buildrel \leftarrow \over x}},
      \,\tau_{n} \le b(n-p k) \right)\,
\end{eqnarray*}

\noindent where $p:=\lfloor {n\over k}\rfloor$. Since $A$ is
bounded, there exists $c_{17}>0$ such that
$\sum_{i=1}^{\nu(y)}\omega(y,y_i)\ge c_{17}$ $\forall y\in \t$. It
yields that
$$
\min_{|x|=p k}P_{\omega}^x\left( \tau_{n}<T_{{\buildrel \leftarrow
\over x}},
      \,\tau_{n} = (n-p k) \right)\ge c_{17}^k\,.
$$

\noindent Hence,
\begin{eqnarray}
     \liminf_{n\rightarrow\infty} {1\over n}\ln\left( P_{\omega}^e\left( \tau_{n}<T_{{\buildrel \leftarrow \over e}},
      \,\tau_{n} \le bn \right)
     \right)
\ge
     -h + {1\over k}\ln(e_k(h,b))\,.\label{LD:infqnched}
\end{eqnarray}

\noindent Take now a general tree $\t$. Notice that since $h\in S$,
$\Q\left( \mathcal{T}_{k,h} \right)>0$ for $k$ large enough, and
there exists almost surely a vertex $z\in \t$ such that the subtree
rooted at it belongs to $\mathcal{T}_{k,h}$. It implies that for
large $k$, (\ref{LD:infqnched}) holds almost surely. Then letting
$k$ tend to infinity and taking the supremum over all $h\in S$ leads
to
\begin{eqnarray*}
     \liminf_{n\rightarrow\infty} {1\over n}\ln\left( P_{\omega}^e\left( \tau_{n}<T_{{\buildrel \leftarrow \over e}},
      \,\tau_{n} \le bn \right)
     \right)
 \ge
 -J_q(b)\,.
\end{eqnarray*}

\noindent For the upper bound in (\ref{LD:iq}), we observe that, for
any integer $k$,
\begin{eqnarray*}
      P_{\omega}^e( \tau_n < T_{{\buildrel \leftarrow \over e}},\, \tau_{n}\le
      bn)
\le
      \sum_{\ell=0}^{ k-1 }
         e^{-\ell n/k}\sum_{ |x| = n }\ind_{A((\ell+1)/k,b,x)}\,.
\end{eqnarray*}

\noindent By Markov's inequality, we have
\begin{eqnarray*}
     \Q\left( \sum_{|x|=n} \ind_{A(h,b,x)} > \left( e(h,b)+ 1/k \right)^n \right)
\le
     { e_n(h,b) \over \left( e(h,b)+ 1/k \right)^n}\le { \left(e(h,b) \over  e(h,b)+ 1/k
     \right)^n}\,,
\end{eqnarray*}

\noindent by (\ref{LD:supehb}). An application of the
Borel--Cantelli lemma proves that $
 \sum_{|x|=n} \ind_{A(h,b,x)} \le \left( e(h,b)+ 1/k \right)^n $ for all but a finite
number of $n$, $\Q$-a.s. In particular, if $e(h,b)+1/k<1$, then $
\sum_{|x|=n} \ind_{A(h,b,x)}= 0$ for $n$ large enough. Consequently,
for $n$ large,
\begin{eqnarray*}
       P_{\omega}^e(\tau_n<T_{{\buildrel \leftarrow \over e}},\, \tau_{n}\le
       bn)
\le
       e^{n/k}k\sup\{e^{-hn}(e(h,b)+1/k)^n,\,h:e(h,b)+1/k\ge
       1\}\,.
\end{eqnarray*}

\noindent We find that
\begin{eqnarray*}
\limsup_{n\rightarrow\infty} {1\over
n}\ln(P_{\omega}^e(\tau_n<T_{{\buildrel \leftarrow \over e}},
\tau_{n}\le bn))\le 1/k + \sup\{-h+ \ln(e(h,b)+1/k),\,
h:e(h,b)+1/k\ge 1\}\,.
\end{eqnarray*}

\noindent Let $k$ tend to infinity and use Corollary \ref{LD:s} to
complete the
proof of (\ref{LD:iq}). \\

\noindent We observe that
\begin{eqnarray*}
    P_{\omega}^e(\tau_n<T_{{\buildrel \leftarrow \over e}}, \tau_{n}\le
    bn)- P_{\omega}^e(\tau_n<T_{{\buildrel \leftarrow \over e}}<\infty,
    \tau_{n}\le bn)
&\le&
    P_{\omega}^e(T_{{\buildrel \leftarrow \over
    e}}=\infty,\,\tau_n\le bn)\\
& \le&
    P_{\omega}^e(\tau_n<T_{{\buildrel \leftarrow \over e}},
    \tau_{n}\le bn)\,.
\end{eqnarray*}

\noindent But $P_{\omega}^e(\tau_n<T_{{\buildrel \leftarrow \over
e}}<\infty, \tau_{n}\le bn)\le P_{\omega}^e(\tau_n<T_{{\buildrel
\leftarrow \over e}}, \tau_{n}\le bn)\max_{i=1,\ldots,\nu(e)}
(1-\beta(e_i))$. Since $\max_{i=1,\ldots,\nu(e)} (1-\beta(e_i))<1$
almost surely, we obtain that
\begin{eqnarray}
\lim_{n\rightarrow\infty} {1\over n}\ln(P_{\omega}^e( \tau_{n}\le
bn)\,|\, D(e)=\infty) = -J_q(b)\,. \label{LD:Jqnchd}
\end{eqnarray}

\noindent In the annealed case, notice that
$\S^e(\tau_n<T_{{\buildrel \leftarrow \over e}}<\infty, \tau_{n}\le
bn)= \S^e(\tau_n<T_{{\buildrel \leftarrow \over e}}, \tau_{n}\le
bn)E_{\P}[1-\beta]$ which leads similarly to
\begin{eqnarray}
\lim_{n\rightarrow\infty} {1\over n}\ln(\S^e( \tau_{n}\le bn)) =
-J_a(b)\,.\label{LD:Janld}
\end{eqnarray}

\noindent We can now finish the proof of the theorem. The continuity
has to be proved only at $b=1$ (since $J_a$ and $J_q$ are convex on
$[1,+\infty[$), which is directly done with the arguments of
\cite{dgpz02} Section 4. We let
$b<1/v=E_{\S^e}[\Gamma_1]/E_{\S^e}[|X_{\Gamma_1}|]$ and we observe
that for any constant $c_{18}>0$,
\begin{eqnarray*}
        \S^e(\tau_n\le b n)
\le
        \S^e( \tau_n<\Gamma_{c_{18}n} ) + \S^e(\Gamma_{c_{18}n} \le b n
        )\,.
\end{eqnarray*}

\noindent Choose $c_{18}$ such that
$b\left(E_{\S^e}[\Gamma_1]\right)^{-1} < c_{18} <
\left(E_{\S^e}[|X_{\Gamma_1}|]\right)^{-1}$. Use Cram\'er's Theorem
with Facts A and B to see that $\S^e( \tau_n<\Gamma_{c_{18}n} )$ and
$\S^e(\Gamma_{c_{18}n} \le b n )$ decrease exponentially. Then,
$\S^e(\tau_n\le b n)$ has an exponential decay and, by
(\ref{LD:Janld}), $I_a(b)>0$ which leads to $I_q(b)>0$ since $I_a\le
I_q$. We deduce in particular that $I_a$ and $I_q$ are strictly
decreasing. Furthermore, $P_{\omega}^e(\tau_{n}\le bn\,|\,
D(e)=\infty)$ tends to $1$ almost surely when $b>1/v$, which in
virtue of (\ref{LD:Jqnchd}), implies that $J_q(b)=0$. By continuity,
$I_q(1/v)=0$ and therefore $I_a(1/v)=0$. Finally, let $a < b$, $
b\in [1, 1/v]$.
\begin{eqnarray*}
P_{\omega}^e\left( an < \tau_n \le bn \,|\, D(e)=\infty \right) =
P_{\omega}^e\left(\tau_n \le bn \,|\, D(e)=\infty
\right)-P_{\omega}^e\left(\tau_n \le an \,|\, D(e)=\infty \right)\,.
\end{eqnarray*}

\noindent Equation (\ref{LD:devupqup2}) follows since $I_q$ is
strictly decreasing. The same argument proves (\ref{LD:devupaup2}).
$\Box$

\subsection{Proof of Theorem \ref{LD:slowdowncond}}

\label{LD:sectslowdowncond}

The proof is the same as before by taking for $b \ge 1$,
\begin{eqnarray*}
\widetilde A
(h,b,x)&:=&\{\omega\,:\,P_{\omega}^e\left(\tau_n=T_x,\,T_{{\buildrel
\leftarrow \over
e}}> \tau_n \ge bn \right)\ge e^{-hn}\}\,,\\
\widetilde {e}_n(h,b)&:=&
E_{\Q}\left[\sum_{|x|=n}\ind_{\widetilde{A}(h,b,x)}\right]\,,
\\
\widetilde{S} &:=& \{h\,:\, \widetilde{e}(h,b)>1\}\,.
\end{eqnarray*}

\noindent Define also for any $b\ge 1$,
\begin{eqnarray*}
       \widetilde{J}_a(b)
&:=&
       -\sup\{-h + \ln({\widetilde e}(h,b))\,,\,h\ge 0\}\,,\\
       \widetilde{J}_q(b)
&:=&
       -\sup\{-h + \ln({\widetilde e}(h,b))\,,\,h\in {\widetilde S}\}\,,
\end{eqnarray*}

\noindent and for any $b\ge 1/v$,
\begin{eqnarray*}
I_a(b) &:=& \widetilde{J}_a(b)\,,\\
I_q(b) &:=& \widetilde{J}_q(b)\,.
\end{eqnarray*}

\noindent We verify that $I_a\le I_q$ and both functions are convex.
We have then for any $b\ge 1$,
\begin{eqnarray}
\lim_{n\rightarrow\infty}{1\over n}\ln\left( \q^e\left(T_{{\buildrel \leftarrow \over e}}> \tau_n \ge b n\right)  \right)
&=& -\widetilde{J}_a(b)\,,  \\
\lim_{n\rightarrow\infty}{1\over n}\ln\left(
P_{\omega}^e\left(T_{{\buildrel \leftarrow \over e}}>  \tau_n \ge b
n\right)  \right) &=& -\widetilde{J}_q(b)  \,.
\end{eqnarray}

\noindent As before, we obtain
\begin{eqnarray*}
\lim_{n\rightarrow\infty}{1\over n}\ln\left( \S^e\left(\tau_n \ge b
n\right)  \right)
&=& -{\widetilde J}_a(b)\,,  \\
\lim_{n\rightarrow\infty}{1\over n}\ln\left( P_{\omega}^e\left(
\tau_n \ge b n\,|\, D(e)=\infty \right)  \right) &=& -{\widetilde
J}_q(b) \,.
\end{eqnarray*}

 \noindent We have ${\widetilde J}_a={\widetilde J}_q=0$ on $[1,1/v]$. In
the case $i>\nu_{min}^{-1}$,  the positivity of $I_a$ and $I_q$ on
$]1/v,+\infty[$ comes from Proposition \ref{LD:timexp} and
Cram\'er's Theorem, which implies that they are strictly increasing.
Equations (\ref{LD:devupaslow2}) and (\ref{LD:devupqslow2}) follow
in that case. In the case $i\le \nu_{min}^{-1}$, we follow the
strategy of \cite{dgpz02}. Let $\eta>0$. As in the proof of
Proposition \ref{LD:timereg}, we set $h_n:=\lfloor
\ln(n)/(6\ln(b))\rfloor$, and for some $b\in \bN$,
\begin{eqnarray*}
w_+  &:=&    \Q\left(\sum_{i=1}^{\nu}A(e_i)\ge {1 + \eta },\,\nu(e) \le b\right)\,, \\
w_-  &:=&    \Q\left(\sum_{i=1}^{\nu}A(e_i)\le {1\over 1+\eta}, \,
\nu(e) \le b\right)\,.
\end{eqnarray*}

\noindent Taking $b$ large enough, we have $w_+>0$ and $w_->0$. We
say that $\t$ is a $n$-good tree if
\begin{itemize}
\item any vertex $x$ of the $h_n$ first generations verifies
 $\nu(x)\le b$ and $\sum_{i=1}^{\nu(x)}A(x_i)\ge {1 +\eta}$\,,
\item any vertex $x$ of the $h_n$ following generations verifies $\nu(x)\le
b$ and $\sum_{i=1}^{\nu(x)}A(x_i)\le {1\over 1 + \eta}$\,.
\end{itemize}
\bigskip
\noindent Then we know that $ Q_n:= \Q(\t\, \mbox{is}\,
n\mbox{-good}) \ge
     \exp(-n^{1/3+o(1)})$. Let $Y'$ be a random walk starting
from zero which increases (resp. decreases) of $1$ with probability
$ {1+\eta \over 2+\eta }$ (resp. ${ 1\over 2+\eta}$). We define
$p'_n$ as the probability that $Y'$ reaches $-1$ before $h_n$. We
show that (\ref{LD:nK}) is still true (by the exactly same
arguments), so that there exists a constant $K>0$ and a
deterministic function $O(n^K)$ bounded by a factor of $n\rightarrow
n^K$, such that
\begin{eqnarray}
      P_{\omega}^e(T_{{\buildrel \leftarrow
      \over e}}>\tau_{2h_n}\ge  n)
\ge
      O(n^K)^{-1}(p'_n)^n\,, \label{LD:goodqnchd}
\end{eqnarray}

\noindent We have, by gambler's ruin formula,
$$
p'_n = 1- {1 \over 1 + \left( {1\over 1+\eta} \right) +\ldots +
\left( {1 \over 1+\eta} \right)^{h_n}} \ge {1 \over 1+\eta}\,.
$$

\noindent Let $k_n:=\lfloor n^{d}\rfloor$ with $d\in (1/3,1/2)$ and
let $f\in (d,1-d)$. We call an $n$-slow tree a tree in which we can
find a vertex $|x| = k_n$ such that $\t_x$ is $n$-good (where $\t_x$
is the subtree rooted at $x$), and for any $y\le x$, we have
$\nu(y)\le \exp(n^f)$. We observe that if a tree is not $n$-slow,
then either there exists a vertex before generation $k_n$ with more
than $\exp(n^f)$ children, or any subtree rooted at generation $k_n$
is not $n$-good. This leads to
\begin{eqnarray*}
\Q(\t\,\mbox{is not}\, n\mbox{-slow}) &\le &
\sum_{\ell=1}^{k_n}E_{GW}[Z_{\ell}]GW(\nu > e^{n^f}) +
E_{GW}\left[ (1-Q_n)^{Z_{k_n}} \right]\\
&\le& k_nm^{k_n}me^{-n^f}+ (1-Q_n)^{(1+\varepsilon)^ {k_n}} +
GW(Z_{k_n}\le (1+\varepsilon) ^{k_n})\,.
\end{eqnarray*}

\noindent We notice that $(1-Q_n)^{(1+\varepsilon)^ {k_n}}\le
\exp(-(1+\varepsilon)^{n^{d+o(1)}})$. Moreover,
\begin{eqnarray*}
GW(Z_{k_n}\le (1+\varepsilon) ^{k_n}) \le
(1+\varepsilon)^{k_n}E_{GW}\left[ {1\over Z_{k_n} }\right]
\end{eqnarray*}

\noindent Observe that for any $k\ge 0$, $E_{GW}\left[ {1\over
Z_{k+1}}\right]\le q_1E_{GW}\left[1\over Z_{k}\right] +
(1-q_1)E_{GW}\left[1\over X_1+X_2\right]$ where $X_1$ and $X_2$ are
independent and distributed as $Z_{k}$. We then verify
$E_{GW}\left[{1\over X_1+X_2}\right]\le {\left(u/2 \right)}\land {
v}$ where $u:= E_{GW}\left[\min(X_1,X_2)^{-1}\right]$ and
$v:=E_{GW}\left[\max(X_1,X_2)^{-1}\right]$. Since $u+v=E_{GW}\left[
{2\over Z_k} \right]$, we deduce that $E_{GW}\left[{1\over
X_1+X_2}\right]\le {2\over 3}E_{GW}\left[ {1\over Z_k} \right]$,
leading to $E_{GW}\left[ {1\over Z_{k+1}}\right]\le (q_1+ {2 \over
3}(1-q_1))E_{GW}\left[ {1\over Z_k} \right]\le (q_1+ {2 \over
3}(1-q_1))^{k+1}$. We get
\begin{eqnarray*}
GW(Z_{k_n}\le (1+\varepsilon) ^{k_n}) \le \left( {(1+\varepsilon )(
q_1+{2 \over 3}(1-q_1))} \right)^{k_n}\,,
\end{eqnarray*}

\noindent and, taking $\varepsilon$ small enough,
\begin{eqnarray}
\Q(\t\,\mbox{is not}\, n\mbox{-slow}) \le \exp(-n^{d+o(1)})\,.
\label{LD:slow}
\end{eqnarray}

\noindent Let $1/v \le a < b$. We want to show that (under the hypothesis $i\le \nu_{min}^{-1}$),
\begin{equation}
\label{LD:qnchzero} \liminf_{n\rightarrow\infty}\ln
P_{\omega}^e({\tau_n\over n}\in [a,b[,\,D(e)>\tau_n)=0\,.
\end{equation}

\noindent If this is proved, the Jensen's inequality gives
\begin{equation}
\nonumber
\liminf_{n\rightarrow\infty}\ln \q^e({\tau_n\over n}\in [a,b[,\,D(e)>\tau_n)=0\,.
\end{equation}

\noindent Equations (\ref{LD:devupqslow2})
and (\ref{LD:devupaslow2}) follow. Therefore, we focus on the proof of (\ref{LD:qnchzero}).\\

\noindent Let $n_1:=n-k_n-2h_n$,  $\delta>0$, and
$G_k:=\{|x|=k\;\mbox{s.t.}\; \t_x\;\mbox{is}\;n\mbox{-slow}\}$. We
have
\begin{eqnarray*}
    \left\{ {\tau_n \over n}\in [a,b[,\,\tau_{{\buildrel
    \leftarrow \over e}}>\tau_n \right\}
\subset
    E_5\cap E_6\cap E_7\,,
\end{eqnarray*}

\noindent with
\begin{eqnarray*}
    E_5 &:=& \left\{{T_{{\buildrel \leftarrow \over e}}}>\tau_{n_1},\, {\tau_{n_1}\over n_1} \in \left[{1 \over v}-\delta, {1\over v}+\delta\right[\right\}\,,\\
    E_6 &:=& \left\{{ X_{\tau_{n_1}}\in
    G_{n_1}}\right\}\,,\\
    E_7 &:=& \left\{ D(X_{\tau_{n_1}})> \tau_n,\,{\tau_n \over n} \in
\left(a-{1 \over v}+\delta,b-{1 \over v}-\delta\right) \right\}\,.
\end{eqnarray*}

\noindent We look at the probability of the event $E_7$ conditioned
on $E_5$ and $E_6$. Therefore, we suppose that $u:=X_{\tau_{n_1}}$
is known, and that the subtree $\t_u$ rooted at $u$ is a $n$-slow
tree. There exists $x_n$ at generation $n_1+k_n$ such that
$\t_{x_n}$ is a $n$-good tree and $\nu(y) \le e^{n^f}$ for any $u\le
y < x_n$. Let also $n$ be large enough so that $k_n \le \delta n$.
It implies that
\begin{eqnarray*}
&& P_{\omega}^u\left( D(u)>\tau_{n},\, {\tau_{n}\over n}\in (a-{1 \over v}+\delta,b-{1
\over v}-\delta) \right) \\
&\ge& P_{\omega}^u\left( D(u)>T_{x_n}=k_n\right)P_{\omega}^{x_n}\left(
D(x_n)>\tau_{n},\, {\tau_{n}\over n}\in (a-{1 \over
v}+\delta,b-{1 \over v}-2\delta) \right)\\
&\ge& \exp(-c_{21}n^{c_{22}})P_{\omega}^{x_n}\left(
D(x_n)>\tau_{n},\, {\tau_{n}\over n}\in (a-{1 \over v}+\delta,b-{1
\over v}-2\delta) \right)  \,, \label{LD:slowqnchd}
\end{eqnarray*}

\noindent for some $c_{22} \in (0,1)$. By definition of a $n$-good
tree, any vertex $x$ descendant of $x_n$ and such that $|x|\le n$
verifies $\nu(x)\le b$. Therefore there exists a constant $c_{23}>0$
such that $P_{\omega}^{y}(\tau_n\le 2h_n)\ge c_{23}^{2h_n}$ for any
$y\ge x_n$, $|y|<n$. By the strong Markov property,
\begin{eqnarray*}
&& P_{\omega}^{x_n}\left( D(x_n)>\tau_{n},\,
{\tau_{n}\over n}\in (a-{1 \over v}+\delta,b-{1 \over
v}-2\delta) \right)\\
 &\ge& P_{\omega}^{x_n}\left(
D(x_n)>\tau_{n},\, {\tau_{n}\over n}\ge a-{1 \over
v}+\delta \right) c_{23}^{2h_n}\,.
\end{eqnarray*}

\noindent Let $L:=a-{1\over v} +\delta$. By equation
(\ref{LD:goodqnchd}),
\begin{eqnarray*}
P_{\omega}^{x_n}\left( D(x_n)>\tau_{n},\,
{\tau_{n}\over n}\ge a-{1 \over v}+\delta \right) \ge O(n^K)^{-1}
\left({1\over 1+\eta}\right)^{L n}\,.
\end{eqnarray*}

\noindent Hence, by the strong Markov property,
\begin{eqnarray*}
\liminf_{n\rightarrow \infty} {1\over n} \ln P_{\omega}^e(E_7\,|\, E_5\cap E_6) &=& \liminf_{n\rightarrow \infty} {1\over n} \ln P_{\omega}^u\left( D(u)>\tau_{n},\, {\tau_{n}\over n}\in (a-{1 \over v}+\delta,b-{1
\over v}-\delta)  \right) \\&\ge& -L(1+\eta)\,.
\end{eqnarray*}

\noindent  This implies that
\begin{eqnarray}
 \lim_{n\rightarrow\infty} {1\over n} \ln P_{\omega}^e\left({\tau_n
\over n}\in [a,b[,\,D(e)>\tau_n\right) \nonumber
&\ge& \liminf_{n\rightarrow} {1\over n} \ln P_{\omega}^e\left(
E_5\cap E_6 \cap E_7
\right) \nonumber \\
&\ge&
 \liminf_{n\rightarrow \infty} {1\over n}  \ln P_{\omega}^e\left( E_5 \cap E_6 \right) - L
 \ln(1+\eta)\,. \label{LD:eqld1}
\end{eqnarray}

\noindent Notice that
\begin{eqnarray*}
E_{\Q}\left[ P_{\omega}^e\left( E_5 \cap E_6^c \right) \right]&=&
E_{\Q}\left[ P_{\omega}^e\left( E_5 \right) -P_{\omega}^e\left( E_5
\cap E_6 \right) \right]\\
&=&\q(E_5)(1-\Q(\t\,\mbox{is}\, n\mbox{-slow}))\\
&\le& \q(E_5)\exp(-n^{d+o(1)})\,,
\end{eqnarray*}

\noindent by equation (\ref{LD:slow}). By Markov's inequality,
$$
\Q(P_{\omega}^e(E_5 \cap E_6^c)\ge {1\over n^2}) \le
n^2\q(E_5)e^{-n^{d+o(1)}}\,.
$$

\noindent The Borel--Cantelli lemma implies that almost surely, for
$n$ large enough,
\begin{eqnarray*}
P_{\omega}^e\left(E_5 \cap E_6\right) \ge P_{\omega}^e(E_5)-{1\over
n^2}\,.
\end{eqnarray*}

\noindent We observe that $P_{\omega}^e(E_5) \rightarrow
P_{\omega}^e(T_{{\buildrel \leftarrow \over e}} =\infty)$ when $n$
goes to infinity. Therefore , equation (\ref{LD:eqld1}) becomes
\begin{eqnarray*}
\lim_{n\rightarrow \infty}{1\over n}\ln\left(
P_{\omega}^e\left({\tau_n \over n}\in
[a,b[,\,D(e)>\tau_n\right)\right)\ge -(a-{1\over v}
+\delta)\ln(1+\eta)\,.
\end{eqnarray*}

\noindent We let $\eta$ go to $0$ to get
\begin{eqnarray*}
\lim_{n\rightarrow \infty}{1\over n}\ln\left(
P_{\omega}^e\left({\tau_n \over n}\in [a,b[,\,D(e)>\tau_n\right)\right) = 0
\end{eqnarray*}

\noindent which proves (\ref{LD:qnchzero}).

\subsection{Proof of Proposition \ref{LD:rootdev}}
\label{LD:sectrootdev}

The speed-up case is quite immediate. Indeed, reasoning on the last
visit to the root, we have
$$
{\q^e( \tau_n \le b n,\,D(e)=\infty)} \le {\q^e( \tau_n \le b n)} \le bn {\q^e(\tau_n \le b n
,\,D(e)=\infty)}\,.
$$

\noindent Therefore, by Theorem \ref{LD:speedupcond},
$$
\lim_{n\rightarrow \infty} {1\over n}\ln \q^e(\tau_n \le b n) =
\lim_{n\rightarrow \infty} {1\over n} \ln \q^e(\tau_n \le b n\,|\,
D(e)=\infty)\,.
$$

\noindent It already gives (\ref{LD:unconspeedupa}) since $I_a$ is
strictly decreasing on $[1,1/v]$. We do exactly the same for the
quenched inequality. Therefore, let us turn to the slowdown case,
beginning with the annealed inequality (\ref{LD:unconslowdowna}). We
follow the arguments of \cite{dgpz02}. We still write $i=\mbox{ess
inf }A$. For technical reasons, we need to distinguish the cases
where $\P(A = i)$ is null or positive. We feel free to deal only
with the case $\P(A = i)=0$, the other one following with nearly any
change. Moreover, we suppose without loss of generality that
$i>\nu_{min}^{-1}$, since the two sides are equal to zero when $i\le
\nu_{min}^{-1}$. Let $k\ge 1$. We write $\ell = k [2]$ to say that
$\ell$ and $k$ have the same parity. Following \cite{dgpz02}, we
write for $b>a>1/v$,
\begin{eqnarray*}
&& P_{\omega}^e(b n>\tau_n \ge a n) \\
&=& \sum_{\ell = k [2]}
\sum_{|x|=k} P_{\omega}^e(b n>\tau_n \ge a n, \tau_n>\ell,\,X_{\ell} =x,
|X_i| >k,
\,\forall \;i=\ell +1\ldots ,\tau_n )\\
&=&\sum_{\ell = k [2]} \sum_{|x|=k} P_{\omega}^e(
\tau_n>\ell,\,X_{\ell} =x)P_{\omega}^{x}(bn-\ell > \tau_n > a
n-\ell,\, D(x)> \tau_n )\,.
\end{eqnarray*}

\noindent By coupling, we have, for $p :=\nu_{min} i >1$,
\begin{eqnarray*}
\sup_{|x|=k} P_{\omega}^e( \tau_n>\ell,\,X_{\ell} =x) \le
P_{\omega}^e(|X_{\ell}| \le k ) \le P( S_{\ell}^{p} \le k )\,,
\end{eqnarray*}

\noindent where $S_{\ell}^p$ stands for a reflected biased random
walk on the half line, which moves of $+1$ with probability $p/1+p$
and of $-1$ with probability $1/1+p$. From (and with the notation
of) Lemma 5.2 of \cite{dgpz02}, we know that for all $\ell$ of the same parity
as $k$,
$$
P(S_{\ell}^p \le k) \le c_k (1 + \delta_k)^{\ell}P(S_{\ell}^p = k,
\, 1\le S_{i}\le k-1,\,i=1,\ldots,\ell-1 )
$$

\noindent where $c_k < \infty$ and $\delta = (\delta_k)$ is a
sequence independent of all the parameters and tending to zero. In
particular, we stress that $\delta$ do not depend on $p$. Hence,
$P_{\omega}^e(b n>\tau_n \ge a n)$ is smaller than
\begin{eqnarray*}
  c_k (1 + \delta_k)^{ bn}\sum_{\ell = k [2]} \sum_{|x|=k}
P(S_{\ell}^p = k, \, 1\le S_{i}\le k-1,\,i=1,\ldots,\ell-1 )W_n(x,\ell)
\end{eqnarray*}

\noindent where
$$
W_n(x,\ell) := P_{\omega}^{x}(bn-\ell
> \tau_n \ge a n-\ell,\, D(x)> \tau_n )\,.
$$

\noindent  We deduce that
\begin{eqnarray}
P_{\omega}^e(b n>\tau_n \ge a n) &\le& c_k (1 + \delta_k)^{b
n}\sum_{\ell = k [2]} \sum_{|x|=k}
P_{\omega_p}^e(\tau_k=\ell,D(e)>\ell)W_n(x,\ell) \nonumber
\\
&=& \nu_{min}^{k}c_k (1 + \delta_k)^{ b n}\sum_{\ell = k [2]}
\sum_{|x|=k} P_{\omega_p}^e(\tau_k=\ell,D(e)>\ell, X_{\ell}=x)
W_n(x,\ell) \,,\label{LD:ldbiased}
\end{eqnarray}

\noindent where $\omega_p$ represents the environment of the
biased random walk on the $\nu_{min}$-ary tree such that for any vertex $x$, $P_{\omega_p}^x(X_1= x_i) ={ p\over {\nu_{min}}( 1 + p)}$ for each child $x_i$, and $P_{\omega}^x(X_1={\buildrel \leftarrow \over x}) = {1 \over 1 +p }$.
Taking the expectations yields that
\begin{eqnarray}
 \q^e(b n>\tau_n \ge a n) \le \nu_{min}^{k}c_k (1 + \delta_k)^{bn}\sum_{\ell = k [2]}
\sum_{|x|=k} P_{\omega_p}^e(\tau_k=\ell,D(e)>\ell, X_{\ell}=x)
E_{\Q}[W_n(x,\ell)]\,.\nonumber \\
\label{LD:rooteq1}
\end{eqnarray}

\noindent Moreover, define
for any $|x|=k$,
$$
\mathcal{S}_{k,\ell}^{+}(\t,x) = \left\{ \{s_i\}_{i=0}^{\ell} :
|s_{i+1}| -|s_i|=1, s_0=0, k-1\ge |s_i| >0, s_{\ell}=x  \right\}
$$

\noindent the set of paths on $\t$ which ends at $x$ in $\ell$ steps
and stays between generation $1$ and $k-1$ before. We notice that,
for any environment $\omega$,
\begin{eqnarray}
\qquad P_{\omega}^e(\tau_{k}=\ell, D(e)>\ell,X_{\ell}=x) =\sum_{
\{s\} \in \mathcal{S}_{k,\ell}^{+}(\t,x)} \sum_{y\in
\t}\omega(y,{\buildrel \leftarrow \over y })^{N(y,{\buildrel
\leftarrow \over y})}\sum_{i=1}^{\nu(y)}\omega(y,y_i)^{N(y,y_i)}
\label{LD:decomp}
\end{eqnarray}

\noindent where for each path $\{s_i\}$, $N(z,y)$ stands for the
number of passage from $z$ to $y$. Let $\varepsilon
>0$, and $\mathcal{G}_k$ denote for any $k$ the set of trees such
that any vertex $x$ of generation less than $k$ verifies
$\nu(x)=\nu_{min}$ and $A(x) \le \mbox{ess inf }A +\varepsilon$. Let
$p':=\nu_{min} (\mbox{ess inf A}+\varepsilon)$. We observe that
\begin{eqnarray*}
P_{\omega_p}^e(\tau_k=\ell,D(e)>\ell, X_{\ell}=x) &=& \sum_{ \{s\}
\in \mathcal{S}_{k,\ell}^{+}(\t,x)} \sum_{y\in \t}\left({1\over
1+p}\right)^{N(y,{\buildrel \leftarrow \over
y})}\sum_{i=1}^{\nu(y)}{\left( p \over \nu_{min}(1+ p)
\right)}^{N(y,y_i)}
\end{eqnarray*}

\noindent Therefore, if $\t$ belongs to $\mathcal{G}_k$, we have by
equation (\ref{LD:decomp}),
\begin{eqnarray*}
P_{\omega_p}^e(\tau_k=\ell,D(e)>\ell, X_{\ell}=x) \le \left({1+
p'\over 1+p}\right)^{\ell}P_{\omega}^e(\tau_{k}=\ell,D(e)>\ell,
X_{\ell}=k)\,.
\end{eqnarray*}

\noindent It entails that
\begin{eqnarray}
&& \ind_{\{ \t \in \mathcal{G}_k \}}\sum_{\ell = k [2]} \sum_{|x|=k}
P_{\omega_p}^e(\tau_k=\ell,D(e)>\ell, X_{\ell}=x)
W_n(x,\ell) \nonumber\\
&\le&
\ind_{\{ \t \in \mathcal{G}_k \}}\left({1+p'\over 1+p}\right)^{b
n}\sum_{\ell = k [2]} \sum_{|x|=k}  P_{\omega}^e( \tau_k
=\ell,D(e)>\ell,\,X_{\ell}=x )W_n(x,\ell) \nonumber  \\
&=&  \ind_{\{ \t \in \mathcal{G}_k \}}\left({1+p'\over 1+p}\right)^{b n}
P_{\omega}^e(b n > \tau_n \ge a n,\,D(e)>\tau_n) \nonumber \\
&\le& \left({1+p'\over 1+p}\right)^{b n} P_{\omega}^e(b n > \tau_n
\ge a n,\,D(e)>\tau_n) \,. \label{LD:ldomegak}
\end{eqnarray}

\noindent Taking expectations gives
\begin{eqnarray}
&& \Q(\t \in \mathcal{G}_k)\sum_{\ell = k [2]} \sum_{|x|=k}
P_{\omega_p}^e(\tau_k=\ell, X_{\ell}=x)
E_{\Q}[W_n(x,\ell)] \nonumber\\
&\le& \left({1+p'\over 1+p}\right)^{b n} \q^e(b n > \tau_n \ge a
n,\,D(e)>\tau_n)\,. \label{LD:eq2part}
\end{eqnarray}

\noindent As before,
\begin{eqnarray*}
&&\q^e(b n > \tau_n \ge a
n,\,D(e)=\infty) + \q^e(b n > \tau_n \ge a
n,\,\infty>D(e)>\tau_n)\\
 &=& \q^e(b n > \tau_n \ge a
n,\,D(e)>\tau_n)\\
&\ge&
\q^e(b n > \tau_n \ge a
n,\,D(e)=\infty)\,.
\end{eqnarray*}

\noindent Since $\q^e(b n > \tau_n \ge a
n,\,\infty>D(e)>\tau_n) \le \q^e(b n > \tau_n \ge a
n,\,D(e)>\tau_n)E_{\Q}[1-\beta]$, we get
$$
\lim_{n\rightarrow \infty} {1\over n} \ln \q^e(b n > \tau_n \ge a
n,\,D(e)>\tau_n) =\lim_{n\rightarrow \infty} {1\over n} \ln \q^e(b n > \tau_n \ge a
n\,|\,D(e)=\infty)\,.
$$

\noindent Consequently, we have by (\ref{LD:rooteq1}) and
(\ref{LD:eq2part})
\begin{eqnarray*}
\limsup_{n\rightarrow\infty}{\q^e(b n>\tau_n\ge an)} \le b\ln\left(
{1+p'\over 1+p} (1+\delta_k) \right) + \lim_{n\rightarrow \infty}
{1\over n} \ln \q^e(b n > \tau_n \ge a n\,|\,D(e)=\infty)\,.
\end{eqnarray*}

\noindent Since $\q^e(cn > \tau_n> bn ) \ge \q^e(c n> \tau_n > b n,\,D(e)=\infty)$, we prove equation (\ref{LD:unconslowdowna}) by taking $p'$ arbitrarily close to $p$, and letting $k$ tend to infinity.\\

\noindent We prove now the quenched equality
(\ref{LD:unconslowdownq}). For any environment $\omega$, construct
the environment $f_p(\omega)$ by setting $A(x) = i $ ($:=$ ess inf
$A$) for any $|x|\le k$. We construct also for $p'>p$, an
environment $f_{p'}(\omega)$ by picking independently $A(x)$ in
$[i,p'/\nu_{min}]$ for any $x \le k$, such that $A(x)$ has the
distribution of $A$ conditioned on $A\in[i,p'/\nu_{min}]$. By
equation (\ref{LD:ldbiased}), we have almost surely
\begin{eqnarray*}
\limsup_{n\rightarrow\infty}{1\over n}  \ln P_{\omega}^e(bn >
\tau_n\ge a n) \le \limsup_{n\rightarrow \infty}{1\over
n}P_{f_p(\omega)}^e(bn > \tau_n\ge an, D(e)>\tau_n) +b
\ln(1+\delta_k)\,.
\end{eqnarray*}

\noindent Equation (\ref{LD:ldomegak}) applied to the environment
$f_{p'}(\omega)$, together with Theorem \ref{LD:slowdowncond} shows
that
\begin{eqnarray*}
\limsup_{n\rightarrow\infty}{1\over n}\ln P_{f_p(\omega)}^e(bn >
\tau_n \ge an, D(e)>\tau_n) \le -I_q(b) + b\ln {1+p' \over 1+p}\,.
\end{eqnarray*}

\noindent Let $p'$ tend to $p$ to get that almost surely,
\begin{eqnarray*}
\limsup_{n\rightarrow\infty}{1\over n} \ln P_{f_p(\omega)}(bn > \tau_n \ge an, D(e)>\tau_n) \le -I_q(b)\,.
\end{eqnarray*}

\noindent Therefore
\begin{eqnarray*}
\limsup_{n\rightarrow\infty} {1\over n} \ln P_{\omega}^e(bn >
\tau_n\ge a n) \le -I_q(b)+ b\ln(1+\delta_k)\,.
\end{eqnarray*}

\noindent When $k$ goes to infinity, we obtain
\begin{eqnarray*}
\limsup_{n\rightarrow\infty} {1\over n} \ln P_{\omega}^e(bn >
\tau_n> an) \le -I_q(b)\,,
\end{eqnarray*}

\noindent which gives equation (\ref{LD:unconslowdownq}).

\subsection{Proof of Proposition \ref{LD:i1}}

\noindent Recall that, for any $\theta \in \r$,
\begin{eqnarray*}
    \psi(\theta)
:=
    \ln\left(E_{\Q}\left[ \sum_{i=1}^{\nu(e)} \omega(e,e_i)^{\theta} \right] \right).
\end{eqnarray*}

\noindent Obviously, for any $n\in \bN$,
$$
{1\over n}\ln\left( \q^e\left(  \tau_n =  n\right)  \right)
=
\ln\left( E_{\Q}\left[ \sum_{i=1}^{\nu(e)} \omega(e,e_i)\right] \right)=\psi(1)\,.
$$

\noindent This proves (\ref{LD:ia1}). For the quenched case, we have
that
$$
P_{\omega}^e\left(  \tau_n =  n\right)
=
\sum_{|x|=n}\prod_{k=0}^{n-1} \omega(x_k,x_{k+1})\,,
$$

\noindent where $x_k$ is the ancestor of the vertex $x$ at
generation $k$. We observe that we are reduced to the study of a
generalized multiplicative cascade, as studied in \cite{Liu00}. The
following lemma is well-known in the case of a regular tree (see
\cite{Fran95} and \cite{vargas}). We extend it easily to a Galton--Watson tree.

\begin{lemma}
We have $\lim_{n\rightarrow \infty}{1\over
n}\ln(\sum_{|x|=n}\prod_{k=0}^{n-1}
\omega(x_k,x_{k+1}))=\inf_{]0,1]}{1\over \theta} \psi(\theta)\,.$
\end{lemma}

\noindent {\it Proof.} When $\psi'(1)<\psi(1)$, Biggins \cite{big77}
shows that $\lim_{n\rightarrow \infty}{1\over
n}\ln(\sum_{|x|=n}\prod_{k=0}^{n-1}
\omega(x_k,x_{k+1}))=\psi(1)=\inf_{]0,1]}{1\over \theta}
\psi(\theta)$. Therefore let us assume that $\psi'(1)\ge \psi(1)$. By the
argument of \cite{Fran95}, we obtain,
\begin{eqnarray*}
\liminf_{n\rightarrow \infty}{1\over
n}\ln\left(\sum_{|x|=n}\prod_{k=0}^{n-1} \omega(x_k,x_{k+1})\right)
\ge \inf_{]0,1]}{1\over \theta} \psi(\theta)\,.
\end{eqnarray*}

\noindent Finally, let $\theta \in ]0,\theta_c[$ where
$\psi(\theta_c)=\inf_{]0,1]}{1\over \theta} \psi(\theta)$. Since
$\left(\sum_{i} a_i \right)^{\theta} \le \sum_{i} a_i^{\theta}$ for
any $(a_i)_i$ with $a_i\ge 0$, it yields that
\begin{eqnarray*}
\limsup_{n\rightarrow \infty}{1\over
n}\ln\left(\sum_{|x|=n}\prod_{k=0}^{n-1} \omega(x_k,x_{k+1})\right)
\le {1\over \theta} \limsup_{n\rightarrow \infty}{1\over
n}\ln\left(\sum_{|x|=n}\prod_{k=0}^{n-1}
\omega(x_k,x_{k+1})^{\theta}\right) \,.
\end{eqnarray*}

\noindent We see that (still by \cite{big77}) $\lim_{n\rightarrow
\infty}{1\over n}\ln(\sum_{|x|=n}\prod_{k=0}^{n-1}
\omega(x_k,x_{k+1})^{\theta}) =\psi(\theta)$. It remains to let
$\theta$ tend to $\theta_c$. $\Box$

\section{The subexponential regime : Theorem \ref{LD:subdev}}

We prove (\ref{LD:subeq}) and (\ref{LD:poleq}) separately. We recall
that the speed $v$ of the walk verifies $v={E_{\S^e}\left[
|X_{\Gamma_1}| \right] \over E_{\S^e}\left[
{\Gamma_1} \right]}$. \\

\noindent {\it Proof of Theorem \ref{LD:subdev} : equation
(\ref{LD:subeq})}. Suppose that either ``$i< \nu_{min}^{- 1} $ and
$q_1=0$'' or ``$i< \nu_{min}^{- 1} $ and  $s < 1$" . Let $a>1/v$ and
$c_{24}>0$ such that $
 c_{24} < \left(E_{\S^e}\left[
X_{\Gamma_1} \right]\right)^{-1} $. We have
\begin{eqnarray*}
      \S^e\left(\tau_n \ge a n\right)
\ge
      \S^e\left(\Gamma_{nc_{24}} \ge a n \right)
      -
      \S^e\left( \Gamma_{nc_{24}} > \tau_n \right)\,.
\end{eqnarray*}

\noindent The second term on the right-hand side decays
exponentially by Cram\'er's Theorem applied to the random walk
$(|X_{\Gamma_n}|,\,n\ge 0)$ (recall that $|X_{\Gamma_1}|$ has
exponential moments by Fact A). The simple inequality
$\S^e\left(\Gamma_{nc_{24}} \ge a n \right)\ge \S^e\left(\Gamma_{1}
\ge a n \right)$ thus implies  by Proposition \ref{LD:timereg} the
lower bound of (\ref{LD:subeq}).
Hence, we turn to the upper bound of (\ref{LD:subeq}). Part (i) of Lemma 6.3 of \cite{dgpz02} states:\\
\noindent {\bf Lemma A~(Dembo et al. \cite{dgpz02}) }{\it Let $Y_1$,
$Y_2,\ldots$ be an i.i.d. sequence with $E(Y_1^2)<\infty$. If
$P(Y_1\ge x)\le \exp(-cx^{\gamma})$ for some $0<\gamma<1$, $c>0$ and
all $x$ large enough, then for all $t>E[Y_1]$,
$$
\limsup_{n\rightarrow\infty}n^{-\gamma}\ln P\left( {1\over n}\sum_{j=1}^n Y_j\ge t \right)\le -c(t-E[Y_1])^{\gamma}\,.
$$
 }
\noindent By Proposition \ref{LD:timereg}, $Y_1=\Gamma_1$ meets the
conditions of the lemma. Therefore, take in lemma A,
$Y_i=\Gamma_{i}-\Gamma_{i-1}$ and $t=a/c_{25}$ where $c_{25}$ is
such that
$$
\left(E_{\S^e}\left[ |X_{\Gamma_1}| \right]\right)^{-1} < c_{25} <a
\left(E_{\S^e}\left[ \Gamma_1 \right]\right)^{-1}.
$$

\noindent In particular, we have $t>E_{\S^e}\left[ \Gamma_1
\right]$. As a result, $\S^e\left(\Gamma_n>tn\right)$ is stretched
exponential. We also know that $\S^e \left(|X_{\Gamma_{nc_{25}}}|
      \le n
      \right)$ is exponentially small by Cram\'er's Theorem ($1/c_{25} < E_{\S^e}\left[
|X_{\Gamma_1}| \right] $). The relation $ \S^e \left( \tau_{n} \ge   a n \right)
\le
      \S^e \left(\Gamma_{nc_{25}} \ge   a
      n \right)
      +
      \S^e \left(|X_{\Gamma_{nc_{25}}}|
      \le n
      \right)$ thus completes the proof. $\Box$\\

\noindent We finish with the case $``\Lambda<\infty"$.\\
\noindent {\it Proof of Theorem \ref{LD:subdev} : equation
(\ref{LD:poleq})}.
 Suppose that $\Lambda<\infty$ and let $a$, $c_{24}$ and $c_{25}$ be as before. We write
\begin{eqnarray*}
\S^e\left(\Gamma_{nc_{24}} \ge a n \right) &\ge&
\sum_{k=1}^{nc_{24}}\S^e\left(\left\{\Gamma_{k}-\Gamma_{k-1} \ge a
n\right\}\cap\left\{ \Gamma_{\ell}-\Gamma_{\ell-1}< an,\,\forall
\ell\neq k \right\}
\right)\\
&=&nc_{24}\S^e\left(\Gamma_1 \ge a n\right)\S^e\left( \Gamma_{1}< an
\right)^{nc_{24}-1}\,.
\end{eqnarray*}

\noindent By Proposition \ref{LD:timepol}, $\S^e\left(\Gamma_1 \ge a
n\right)=n^{-\Lambda +o(1)}$. Therefore $\S^e\left( \Gamma_{1}< an
\right)^{nc_{24}-1}$ tends to $1$ (since $\Lambda>1$). Consequently,
\begin{eqnarray*}
\S^e\left(\Gamma_{nc_{24}} \ge a n \right) \ge n^{1-\Lambda+o(1)}\,,
\end{eqnarray*}

\noindent which gives the lower bound of (\ref{LD:poleq}), by the
inequality $      \S^e\left(\tau_n \ge a n\right) \ge
      \S^e\left(\Gamma_{nc_{24}} \ge a n \right)
      -
      \S^e\left( \Gamma_{nc_{24}} > \tau_n \right)$. Turning, to the upper bound, write as before
$ \S^e \left( \tau_{n} \ge   a n \right)
\le
      \S^e \left(\Gamma_{nc_{25}} \ge   a
      n \right)
      +
      \S^e \left(|X_{\Gamma_{nc_{25}}}|
      \le n
      \right)$. We already know that $\S^e \left(|X_{\Gamma_{nc_{25}}}|
      \le n
      \right)$ is exponentially small. Let $H_n:=\Gamma_n-E_{\S^e}[\Gamma_1]n$. When $E[H_1^p]<\infty$,
example 2.6.5 of \cite{Petrov75} says that if $p\ge 2$,
$$
P(H_n>x)\le \left(1+2/p\right)^pnE[H_1^p]x^{-p}+
\exp(-2(p+2)^{-2}e^{-p}x^2/(nE[H_1^2]))
$$
and example 2.6.20 of \cite{Petrov75}, combined with Chebyshev's
inequality, shows that if $1\le p\le 2$,
$$
P(H_n>x)\le \left(2-1/n\right)nE[H_1^p]x^{-p}\,.
$$

\noindent By Proposition \ref{LD:timepol}, $E[H_1^p]<\infty$, for
any $p<\Lambda$. We take $x=({a\over
c_{25}}E_{\S^e}[|X_{\Gamma_1}|]-E_{\S^e}[\Gamma_1])n$ to see that $
\S^e(\Gamma_{nc_{25}} \ge an) \le c(p)n^{1-p} $ for any $p<\Lambda$.
Let $p$ tend to $\Lambda$ in order to complete the proof of equation
(\ref{LD:poleq}). $\Box$

\bigskip

\noindent {\bf Acknowledgements}: I am very grateful of Zhan Shi for
 his guidance in the redaction of this work. I thank also the anonymous referee
for suggesting me several improvements.

\bibliographystyle{plain}
\bibliography{biblio}

\bigskip


{\footnotesize

\baselineskip=12pt

\hskip110pt Laboratoire de Probabilit\'es et Mod\`eles Al\'eatoires

\hskip110pt Universit\'e Paris VI

\hskip110pt 4 Place Jussieu

\hskip110pt F-75252 Paris Cedex 05

\hskip110pt France

\hskip110pt {\tt aidekon@ccr.jussieu.fr}

}

\end{document}